\documentclass[12pt]{article}
\usepackage[utf8]{inputenc}
\usepackage{extsizes}
\usepackage{amsmath}
\usepackage{environ}
\usepackage[left=0.65in,right=0.65in,top = 0.7in, bottom = 0.7in]{geometry}
\usepackage{graphicx}
\usepackage[toc,page]{appendix}
\usepackage[inline]{enumitem}
\usepackage{multicol}
\usepackage[T1]{fontenc}
\usepackage{multirow}
\usepackage{amssymb}
\usepackage{setspace}
\usepackage{titling}
\usepackage{bigints}
\usepackage{bibentry}
\usepackage{import}
\usepackage{subcaption}
\usepackage{upgreek}
\usepackage{longtable}
\usepackage{bm}
\usepackage{bbm}
\usepackage{algorithm}
\usepackage{siunitx}
\usepackage{algpseudocode}
\usepackage{mathtools, nccmath}
\usepackage[mathlines]{lineno}
\usepackage[most]{tcolorbox}
\usepackage{tabularray}
\usepackage{esint}
\usepackage[unicode=true,hypertexnames=false]{hyperref}
\usepackage[nameinlink, capitalise]{cleveref}
\usepackage{nameref}
\usepackage{booktabs}
\usepackage{collcell}
\usepackage{autobreak}
\usepackage{listings}
\usepackage{pgf,tikz,pgfplots}
\usepackage{asymptote}
\usepackage{fancyhdr}
\usepackage{lscape}
\usepackage{authblk}
\usepackage{pdfpages}
\usepackage{soul}
\pgfplotsset{compat=1.15}
\usepackage{amsthm}
\usepackage{mathrsfs}
\numberwithin{equation}{section}
\usepackage{float}
\usepackage{longfigure}
\usepackage{adjustbox}
\usepackage[english]{babel}
\usetikzlibrary{arrows}
\usetikzlibrary{shapes}
\makeatletter
\newcommand*{\currentname}{\@currentlabelname}
\usepackage{titlesec}
\titleclass{\subsubsubsection}{straight}[\subsection]
\newcounter{subsubsubsection}[subsubsection]
\renewcommand\thesubsubsubsection{\thesubsubsection.\arabic{subsubsubsection}}

\titleformat{\subsubsubsection}
  {\normalfont\normalsize\bfseries}{\thesubsubsubsection}{1em}{}
\titlespacing*{\subsubsubsection}
{0pt}{3.25ex plus 1ex minus .2ex}{1.5ex plus .2ex}
\def\toclevel@subsubsubsection{4}
\def\l@subsubsubsection{\@dottedtocline{4}{7em}{4em}}
\makeatother

\usepackage{upquote}
\usepackage{chngcntr}
\usepackage{apptools}
\usepackage{etoolbox}
\AtAppendix{\counterwithin{theorem}{section}}

\newcommand\mydots{\ifmmode\ldots\else\makebox[1em][c]{.\hfil.\hfil.}\thinspace\fi}

\usepackage{scalerel,stackengine}
\stackMath
\definecolor{qqwuqq}{rgb}{0.,0.39215686274509803,0.}
\definecolor{ududff}{rgb}{0.30196078431372547,0.30196078431372547,1.}
\newsavebox{\measurebox}
\newcommand{\fitRHS}[1]{%
    \sbox{\measurebox}{$\displaystyle #1$}%
    \ifdim\wd\measurebox > 0.9 \linewidth
        \resizebox{0.9 \linewidth}{!}{\usebox{\measurebox}}%
    \else
        \usebox{\measurebox}%
    \fi
}
\usepackage{calc}
\newsavebox{\resizeboxaltbox}
\newcommand{\resizeboxalt}[3]{%
\sbox{\resizeboxaltbox}{#3}%
\ifdim\wd\resizeboxaltbox>#1\relax
    \resizebox{#1}{#2}{\usebox{\resizeboxaltbox}}%
\else
    \usebox{\resizeboxaltbox}%
\fi
}
\counterwithin{algorithm}{section}
\newcommand{\ignore}[1]{}
\newcommand{\nobibentry}[1]{{\let\nocite\ignore\bibentry{#1}}}

\newtcbtheorem[auto counter]{theo}%
  {Theorem}{fonttitle=\bfseries\upshape,
     arc=0mm, colback=black!5!white,colframe=black!75!black}{theorem}
\hypersetup{
    pdftitle={Algorithmic Non-dimensionalisation of Polynomial Dynamical Systems with CRNNondim},
    pdfauthor={Georgio Hawi and Martin Wechselberger},
    colorlinks=true,
    linkcolor=blue, 
    citecolor=blue,
    filecolor=magenta,
    urlcolor=cyan
    }

\date{}
\title{Algorithmic Non-dimensionalisation of Polynomial Dynamical Systems with CRNNondim}
\author[1, *]{Georgio Hawi}
\author[1]{Martin Wechselberger}
\affil[1]{School of Mathematics and Statistics, University of Sydney, Sydney, Australia}
\affil[*]{Corresponding author: \href{mailto:georgio.hawi@sydney.edu.au}{georgio.hawi@sydney.edu.au}}
\floatstyle{plaintop}
\restylefloat{table}
\newtheorem{theorem}{Theorem}[section]
\newtheorem{lemma}{Lemma}[section]
\newtheorem{corollary}{Corollary}[section]
\newtheorem{claim}{Claim}
\newtheorem{proposition}{Proposition}[section]

\theoremstyle{definition}

\newtheorem{numberedremark}{Remark}[section]
\crefname{numberedremark}{remark}{remarks}
\Crefname{numberedremark}{Remark}{Remarks}
\newtheorem{definition}{Definition}

\newtheorem{example}{Example}[section]

\newcommand{\set}[1]{\left \{ #1 \right \}}

\DeclareMathOperator{\sign}{sgn}

\allowdisplaybreaks
\newbool{forpublication}
\usepackage[numbers,square,sort&compress]{natbib}
\usepackage{bibunits}
\usepackage{hypernat}

\makeatletter
\pretocmd{\@startbibunit}{%
  \begingroup
    \count@\@bibunitauxcnt
    \advance\count@\@ne
    \xdef\@extra@binfo{.bu\the\count@}%
    \xdef\@extra@b@citeb{.bu\the\count@}%
  \endgroup
}{}{}
\apptocmd{\endbibunit}{%
  \gdef\@extra@binfo{}%
  \gdef\@extra@b@citeb{}%
}{}{}
\makeatother
\begin{document}
\maketitle
\begin{abstract}
Although non-dimensionalisation is a fundamental step in modelling, it is rarely obvious, especially for large dynamical systems, how to systematically construct optimal non-dimensionalisations. In this work, we develop a general mathematical framework for non-dimensionalising polynomial ODEs by interpreting them as chemical reaction networks (CRNs) and characterise the space of admissible scalings. Within this viewpoint, admissible scalings of time and species are encoded by Laurent monomials, and we derive exact formulae for the minimal number of multiplicative dimensionless parameters required by the reaction-channel coefficients and by the full initial value problem, together with characterisations of the scale assignments that attain structural lower bounds. These counts are expressed in terms of structural properties of the underlying CRN, and we show that the full initial value problem requires no more parameters than the corresponding Buckingham count. These results lead to constructive algorithms which we analyse rigorously and implement in the open-source Python package CRNNondim, available at \url{https://github.com/georgio-hawi/crnnondim}. Given a polynomial system or CRN, CRNNondim determines species and time scalings that bring variables as close to unity as possible while minimising the number of free parameters, and among the admissible scalings it chooses free parameters that result in the corresponding dimensionless system being sparse and simple. We apply CRNNondim to numerous well-known biological and biochemical models, yielding simpler dimensionless systems with fewer free parameters, and interpreting the parameters that govern their dynamics.
\end{abstract}
\begin{bibunit}[vancouver]
\section{Introduction}
Polynomial ordinary differential equation (ODE) systems form one of the most ubiquitous modelling frameworks in applied mathematics and the sciences, with applications in biology, physics, and chemistry \cite{Strogatz2018, Hethcote2000, Lorenz2004, Volterra1926, Horn1972, kermack1927contribution, Nagumo1962}. In their simplest form, such systems are expressible as
\begin{equation}
  \frac{dX_i}{dT} = f_i(\mathbf{X};\theta), \qquad i = 1,\dots,n, \label{eq:poly-ode}
\end{equation}
where each component $f_i$ is a polynomial in the state variables $\mathbf{X}=(X_1,\dots,X_n)\in\mathbb{R}^n$, $T$ represents time, and $\theta$ denotes positive physical parameters such as rate constants. \\~\\
Polynomial systems such as \eqref{eq:poly-ode} are assumed to be dimensionally consistent and are typically formulated in dimensional variables. For example, the state vector $\mathbf{X}$ carries units (such as concentrations, densities, or population sizes), the time variable $T$ is measured in physical units (such as seconds, hours, or days), and the parameters $\theta$ possess heterogeneous physical dimensions. Even for moderately complex systems, this dimensional formulation obscures the effective number of degrees of freedom and can hide the relative importance of processes operating on disparate time and concentration scales. The canonical approach to address this dimensional heterogeneity is via non-dimensionalisation, so that by introducing characteristic scales for state variables and parameters, \eqref{eq:poly-ode} is reformulated in terms of dimensionless quantities.\\~\\
The theoretical foundation for non-dimensionalisation is established by the classical Buckingham $\Pi$ theorem \cite{Buckingham1914}, which asserts that given a list of $n$ physical quantities composed of $k$ fundamental dimensions (such as mass, length, and time), any dimensionally homogeneous relation among them can be expressed in terms of $n-m$ independent dimensionless combinations, known as $\Pi$-groups, where $m$ is the rank of the dimensional matrix. Modern treatments formalise this by constructing a mapping from variables to a dimensional matrix, so that by identifying a basis for the null space of this matrix, one obtains the specific exponent vectors required to form the dimensionless products \cite{Barenblatt2003, Curtis1982}. In the context of dynamical systems, this construction results in dimensionless variables $\mathbf{x}$, parameters $\tilde{\theta}$, and time $t$ such that \eqref{eq:poly-ode} can be rewritten in the dimensionless form
\begin{equation}
  \frac{dx_i}{dt} = \tilde{f}_i(\mathbf{x};\tilde{\theta}), \qquad i=1,\dots,n, \label{eq:dimensionless-system}
\end{equation}
where $\tilde{f}_i$ depends only on dimensionless quantities. The mapping from $(\mathbf{X},T,\theta)$ to $(\mathbf{x},t,\tilde{\theta})$ is typically constructed by introducing characteristic scales for $\mathbf{X}$ and $T$ and scaling accordingly \cite{Barenblatt2003}.\\~\\
Importantly, non-dimensionalisation is more than a cosmetic change of variables; it is a foundational step in the qualitative and quantitative analysis of complex systems. Non-dimensionalisation can reduce the system to a smaller set of independent parameter combinations, revealing the true degrees of freedom and isolating dimensionless ratios that control qualitative behaviour. In particular, the identification of naturally small or large dimensionless groups allows for perturbation analysis and rigorous model reductions, which can be formalised within the framework of geometric singular perturbation theory (GSPT); see, for example, \cite{Feliu2022,Goeke2015,Goeke2013,Lapuz2026} for more details. Additionally, by rescaling concentrations and time so that typical dimensionless trajectories have comparable magnitudes, the process can mitigate ill-conditioning and enhance the stability and accuracy of numerical integration, parameter estimation, optimisation, and uncertainty quantification algorithms \cite{Hairer1996}.\\~\\
Despite its conceptual simplicity, non-dimensionalisation is rarely carried out in a systematic, algorithmic manner for large nonlinear models. Classical expositions of the Buckingham $\Pi$ theorem emphasise small illustrative examples, often drawn from continuum mechanics and fluid dynamics \cite{bridgman1922dimensional}, while in mathematical biology and chemical kinetics, standard references typically rely on physical intuition to manually select scales for a handful of paradigmatic models \cite{murray2002mathematical, Lin1988}. Consequently, the construction of $\Pi$-groups is often performed by inspection---a manual approach that assumes a small number of variables and is error-prone, difficult to reproduce, and ill-suited for the large systems routinely encountered in modern modelling \cite{szirtes2007applied, Hart1995}. From a structural perspective, standard treatments focus on variable mapping but treat the equations as a ``black box'', failing to exploit the specific algebraic structure of the underlying vector fields. Moreover, the choice of characteristic scales is rarely unique and may rely on subjective, context-specific intuition that lacks rigorous definition.\\~\\
This being said, greater attention to algorithmic approaches to non-dimensionalisation and scaling reduction has been given in recent years. Hubert and Sedoglavic computed reductions from dilatation and translation Lie-point symmetries with polynomial arithmetic complexity \cite{HubertSedoglavic2006}, and Hubert and Labahn formulated scaling invariants and symmetry reduction using integer linear algebra \cite{HubertLabahn2013}. Toric-ideal methods provide another algebraic route to primitive dimensional invariants \cite{Atherton2014}, and Tanburn et al. developed an automated method for rational first-order ODE systems based on scaling-symmetry reduction \cite{Tanburn2025}. In this work, we approach automating non-dimensionalisation by developing a rigorous mathematical framework, algorithms and software that exploit the special structure of polynomial ODE systems by representing them as chemical reaction networks (CRNs). We also consider the complementary problem of choosing species scales from prescribed, physically meaningful reference quantities, typically the nonzero initial concentrations, and determining how many multiplicative dimensionless parameters are then required.\\~\\
A CRN is fundamentally defined by a finite set of species, a finite set of complexes (each formally a linear combination of species with non-negative integer coefficients), and a finite set of reactions between complexes. In particular, the network structure determines the stoichiometry, and we assume that the dynamical behaviour of the concentration vector $\mathbf{X}(T)\in \mathbb{R}^n$ abides by mass-action-like kinetics
\begin{equation}
  \frac{d\mathbf{X}}{dT} = SV(\mathbf{X}), \label{ministoich}
\end{equation}
where $S\in\mathrm{Mat}_{n,r}(\mathbb{Z})$ is the stoichiometric matrix, whose $j$-th column represents the net stoichiometric change associated with reaction $j$, and $V(\mathbf{X})=(V_1(\mathbf{X}),\dots,V_r(\mathbf{X}))^\intercal$ is the vector of reaction rates. Every polynomial vector field can be formally decomposed into reaction-like monomial terms. As such, \eqref{eq:poly-ode} can always be algebraically expressed in the form of \eqref{ministoich}, with each rate taking the form of a monomial $V_j(\mathbf{X}) = k_j \prod\limits_{i=1}^n X_i^{\nu_{ij}}$, where $k_j>0$ is the rate constant and $\nu=(\nu_{ij})$ is a reactant matrix \cite{Horn1972}.\\~\\
However, the physical validity of such a network depends on the signs of the terms in \eqref{eq:poly-ode}. A polynomial system is known as kinetic if each $f_i(\mathbf{X},\theta)$ can be written as a sum of monomials such that whenever the coefficient of $\mu$ in $f_i$ is negative, the monomial $\mu$ contains $X_i$ as a factor with a positive exponent, ensuring that species are only consumed when they are present. The ``Hungarian lemma'' establishes that while any polynomial system has a formal reaction-like representation, only such kinetic systems can be realised as a mass-action CRN \cite{hars1981inverse}. Conversely, any polynomial system satisfying the condition admits at least one mass-action realisation \cite{hars1981inverse, Craciun2020}. These systems inherit the powerful structural theory pioneered by Horn, Jackson, and Feinberg, known as chemical reaction network theory (CRNT), which rigorously connects graph-theoretic properties of mass-action CRNs to their structural properties, including multistationarity, stability, reversibility and persistence \cite{Horn1972, Feinberg1987, Feinberg1988, Feinberg2019, Craciun2013}.\\~\\
From the perspective of non-dimensionalisation, viewing polynomial systems through the lens of CRNs is particularly advantageous because the nonzero entries of $S$ determine where each reaction monomial contributes, while the monomial structure of $V(\mathbf{X})$ records which species scales multiply each rate constant. These interact with the dimensional information carried by species and rate constants in a highly structured way that classical dimensional analysis does not fully exploit. By working with CRNs instead of the polynomial systems directly, non-dimensionalisation can be formulated as a constrained integer linear algebra problem on the exponents of these monomials. Importantly, the algebraic results use the fixed matrices $S$ and $\nu$ and the assignment of rate symbols to reaction channels, and therefore also apply to formal reaction-like decompositions that do not satisfy the kinetic condition. \\~\\
In this work, we develop a rigorous mathematical framework for non-dimensionalising CRNs and, by extension, general polynomial dynamical systems, providing constructive algorithms for non-dimensionalising such systems, and incorporating these into the open-source Python package CRNNondim. For a fixed reaction-channel decomposition, these algorithms first minimise the number of multiplicative parameters and then select characteristic scales that keep the dimensionless variables as close to unity as possible over a specified time interval. The resulting parameter bases are simplified by unimodular transformations, aiming to obtain a dimensionless system that is sparse and simple.
\section{Setup}\label{setupsection}
For convenience, we provide a collection of principal notation used in the remainder of this work in Table A.1. Throughout this work, we assume that $0$ is included in the set of natural numbers. 
\subsection{Law of Mass Action}
We consider a CRN with ODEs given by
\begin{equation}
\frac{d\mathbf{X}}{dT} = SV(\mathbf{X}), \label{maincrnform}
\end{equation}
where $\mathbf{X}\in\mathbb{R}^n_{\geq 0}$ represents a vector of dimensional chemical species concentrations involved and $T$ is dimensional time. Additionally, $S$ is referred to as the \textbf{stoichiometric matrix}---an $n \times r$ matrix of integers encoding the reaction network structure, where $r$ is the number of reactions. Finally, $V(\mathbf{X}) \in \mathbb{R}^{r}_{\geq 0}$ is the corresponding reaction-rate vector which encodes relevant reactions alongside their rate constants, which, by the law of mass action, can be expressed in the form
\begin{equation}
V_j(\mathbf{X}) = k_j \prod_{i=1}^n X_i^{\nu_{ij}},\quad k_j>0, \nu_{ij}\in \mathbb{N}, \label{eq:rates}
\end{equation}
where $k_j$ is the rate constant corresponding to reaction $j$, and $\nu=(\nu_{ij})$ is a reactant matrix \cite{Horn1972}. For a physical mass-action CRN realisation, we additionally require that $\nu_{\cdot j}+S_{\cdot j}\in\mathbb{N}^n$ for every reaction $j$, so that each product complex is valid, though this additional condition is not required for the algebraic results henceforth. In addition, multiple reactions may share the same rate symbol, and we denote by $\bm{\kappa} := \{\kappa_1,\dots,\kappa_{h}\}$ the set of distinct rate symbols among the reactions, i.e. there exists a surjective map $\rho:\set{1,\dots,r} \to \set{1,\dots,h}$ such that $k_j = \kappa_{\rho(j)}$. Throughout this work, the CRN realisation and reaction-channel decomposition are fixed, and all rate and scale symbols are treated as algebraically independent formal indeterminates, so that there are no relationships between distinct rate constants. Furthermore, we denote $P=|\sign(S)| \in \mathrm{Mat}_{n,r}(\set{0,1})$ as the \textbf{active participation matrix}, so that $P_{ij}=1$ if and only if $S_{ij}\neq 0$. In particular, we call $(i,j)$ an \textbf{active species-reaction pair} when $P_{ij}=1$, and denote the set of such pairs by $\mathcal{P}$. \\~\\
We assume that $n,r\geq 1$ and that each column of $S$ is nonzero, so that every reaction has at least one active species, thereby allowing zeroth-order reactions of the form $\emptyset \to X_i$. We also assume that every species is active in at least one reaction, so that each row of $P$ is nonzero; equivalently, for each $i$ there exists $j$ with $S_{ij}\neq 0$. We note that a species that occurs in the rate monomial with $\nu_{ij}>0$ but has $S_{ij}=0$, as may occur for a catalyst, is not active in that reaction; nevertheless, its scale still enters the coefficients through the factor $\prod_{\ell=1}^n K_\ell^{\nu_{\ell j}}$.
\subsection{Non-dimensionalisation}
We can non-dimensionalise \eqref{maincrnform} by rescaling time and species via
\begin{align}
    T &= K_{T}t, \label{Tscaling}\\
    X_i &= K_i x_{i}, \quad \text{for $i \in \mathcal{I}:=\set{1,\dots,n}$,} \label{Xscaling}
\end{align}
where $K_i>0$ for all $i \in \mathcal{I}$ are \textbf{species scales}, $K_T>0$ is a time scale, and we denote $\mathbf{x}=[x_1,x_2,\dots,x_{n}]^\intercal$ \cite{Buckingham1914, Strogatz2018}. Importantly, $t$ and all elements in $\mathbf{x}$ are dimensionless. \\~\\
Writing $\mathcal{K}=\operatorname{diag}((K_i)_{i\in\mathcal{I}})$, positivity gives $\mathbf{X}=\mathcal{K}\mathbf{x}$ with $\mathcal{K}^{-1}=\operatorname{diag}((K_i^{-1})_{i\in\mathcal{I}})$, and hence
\begin{equation}
  \frac{d\mathbf{x}}{dt} = \mathcal{K}^{-1}\frac{d\mathbf{X}}{dt} = \mathcal{K}^{-1}K_T\frac{d\mathbf{X}}{dT} = \mathcal{K}^{-1}SK_TV(\mathcal{K}\mathbf{x}) = \tilde{S}\tilde{V}(\mathbf{x}). \label{nondimform}
\end{equation}
Substituting \eqref{eq:rates} into \eqref{nondimform} reveals that for each active species-reaction pair $(i,j)$, the contribution of reaction $j$ to the coefficient of the monomial $\mathbf{x}^{\nu_{\cdot j}}:=\prod_{m=1}^n x_m^{\nu_{mj}}$ in the $i\textsuperscript{th}$ component of \eqref{nondimform} is
\begin{equation}
    \hat{c}_{ij} = \frac{K_T \kappa_{\rho(j)}}{K_i} \left(\prod_{\ell=1}^n K_\ell^{\nu_{\ell j}}\right)S_{ij}. \label{fullcijeqn}
\end{equation}
In particular, this implies that the set of all nonzero such monomial contributions, ignoring integer factors due to $S_{ij}$, is given by 
\begin{equation}
  \set{{\tilde{c}}_{ij}} = \left\{\frac{K_T \kappa_{\rho(j)}}{K_i} \left(\prod_{\ell=1}^n K_\ell^{\nu_{\ell j}}\right)\middle|P_{ij}=1\right\}. \label{KTsetcijeqn}
\end{equation}
Since a factor of $K_T$ occurs for all coefficients in $\{\tilde{c}_{ij}\}$, we suppress this common factor for the structural analysis and reintroduce it below. Thus, we denote
\begin{equation}
  \set{{c}_{ij}} = \left\{\frac{\kappa_{\rho(j)}}{K_i} \left(\prod_{\ell=1}^n K_\ell^{\nu_{\ell j}}\right)\middle|P_{ij}=1\right\}. \label{setcijeqn}
\end{equation}
Let $\{a_q\}_{q\in U_\mathrm{cand}}$ be a finite set of positive \textbf{scale symbols}, typically the nonzero initial concentrations, possibly augmented by other reference quantities. We assume that there exists a map, called a \textbf{scale assignment}, $\sigma:\mathcal{I}\to U_\mathrm{cand}$ where $K_i=a_{\sigma(i)}$. We write $U:=\sigma(\mathcal{I})$ for the set of used scale indices and $s:=|U|$, identifying $\sigma$ with its surjective corestriction to $U$ when convenient. Let $\mathcal{J}:=\{i\in\mathcal{I}\mid X_i(0)>0\}$, and assume that $X_i(0)=a_{\iota(i)}$ for an initial-value map $\iota:\mathcal{J}\to U_\mathrm{cand}$. The dimensionless initial conditions are then
\begin{equation}
 x_i(0)=\frac{a_{\iota(i)}}{a_{\sigma(i)}}\quad(i\in\mathcal{J}),\qquad x_i(0)=0\quad(i\notin\mathcal{J}). \label{eq:dimlessIC}
\end{equation}
In practice, it is natural to set $K_i=X_i(0)$ if $X_i(0)>0$, so that $x_i(0)=1$, and to choose the remaining rescaling factors from this set of $K_i$ whenever dimensionally admissible. This normalises the nonzero initial values, although the magnitudes attained later depend on the dynamics.\\~\\
Let $\mathbb{Z}^k$ denote the group of exponent vectors with respect to $k$ fundamental dimensions, one of which is time, with basis vector $e_T$. Each species and each scale symbol carries a dimension $\dim(X_i),\dim(a_q)\in\mathbb{Z}^k$ with zero time component. We assume dimensional consistency, so that every rate symbol $\kappa_p$ has a dimension satisfying
\begin{equation}
\dim(\kappa_{\rho(j)})+\sum_{\ell=1}^n\nu_{\ell j}\dim(X_\ell)-\dim(X_i)=-e_T\qquad\text{for all }(i,j)\in\mathcal{P}, \label{eq:dimconsistency}
\end{equation}
and $\dim(a_{\iota(i)})=\dim(X_i)$ for $i\in\mathcal{J}$. A scale assignment is \textbf{admissible} if $\dim(a_{\sigma(i)})=\dim(X_i)$ for every $i$, and an \textbf{admissible family} $\Sigma$ is a nonempty set of such assignments, possibly subject to additional user-specified restrictions. When these restrictions are imposed separately on species, we write $\sigma(i)\in\mathcal{A}_i\subseteq U_\mathrm{cand}$, with $\mathcal{A}_i$ containing only dimensionally compatible symbols. Furthermore, a time scale is a positive Laurent monomial in the rate and scale symbols with the dimension of time, with all $\tilde{c}_{ij}$ and all nonzero values in \eqref{eq:dimlessIC} being dimensionless for admissible species and time scales. Moreover, subtracting \eqref{eq:dimconsistency} for two active species of the same reaction shows that they have the same dimension, so the species dimension is constant on every connected component of $G$.
\subsection{The Species Co-participation Graph, Scale Graph, and Stoichiometric Rate Lattice}
We define the \textbf{species co-participation graph}, $G$, to be the simple undirected graph on vertices $\{1,\dots,n\}$ with an edge $\set{i, i^\prime}$ if there exists a reaction $j$ such that $P_{ij} = P_{i^\prime j}=1$. Let $C_1, \dots, C_q$ be the connected components of $G$, so that $\bigcup_{\alpha=1}^q C_\alpha = \mathcal{I}$, and for each component $C_\alpha$, we define the \textbf{used scale index set} $U_\alpha:=\set{\sigma(i)\mid i\in C_\alpha}\subseteq U$.
We also define the \textbf{scale graph}, $H_\mathrm{scale}$, as the undirected graph whose vertex set is $U$, formed by the union of the cliques on each $U_\alpha$, so that, with a slight abuse of notation,
\begin{equation}
    H_\mathrm{scale} := \bigcup_{\alpha=1}^q K_{U_\alpha}, \label{hscaledef}
\end{equation}
where $K_{U_\alpha}$ denotes the clique on vertex set $U_\alpha$. Intuitively, two scale symbols are connected if they co-occur somewhere within the same species component of $G$. Additionally, we denote $c_\mathrm{scale}$ as the number of connected components of $H_\mathrm{scale}$, and let $\mathcal{K}_{\mathrm{comp}} = \set{1,\dots,c_\mathrm{scale}}$.\\~\\
We define $\psi: \mathcal{I} \to \mathcal{K}_{\text{comp}}$ to be the component map such that $\psi(i)$ is the index of the connected component of $H_{\mathrm{scale}}$ containing the scale symbol $\sigma(i)$. We also define the map $\psi_V: \{1,\dots,r\} \to \mathcal{K}_{\text{comp}}$ such that $\psi_V(j)$ is the component index of any active species in reaction $j\in \set{1,\dots,r}$, which is well-defined since all active species in a given reaction belong to the same connected component of $G$. Accordingly, we define the \textbf{reaction scale vector} $z_j$ for reaction $j$ as
\begin{equation}
    z_j := \left( \sum_{\ell=1}^n \nu_{\ell j} w_{\psi(\ell)} \right) - w_{\psi_V(j)}  \in \mathbb{Z}^{c_{\mathrm{scale}}},
\end{equation}
where $w_p \in \mathbb{Z}^{c_{\mathrm{scale}}}$ denotes the $p$-th standard basis vector of $\mathbb{Z}^{c_{\mathrm{scale}}}$. We note that the sum in $z_j$ uses the full reactant exponent vector $\nu_{\cdot j}$, so species with $\nu_{\ell j}>0$ and $S_{\ell j}=0$ still contribute through their scale components. Then, the \textbf{stoichiometric rate lattice}, $\Lambda_{\mathrm{rate}} \subseteq \mathbb{Z}^{c_{\mathrm{scale}}}$, is the $\mathbb{Z}$-lattice generated by the differences of these $z_j$ for reactions sharing rate constants, so that
\begin{equation}
    \Lambda_{\mathrm{rate}} := \operatorname{span}_{\mathbb{Z}} \left( z_j - z_k \middle| j,k \in \{1,\dots,r\} \text{ such that } \rho(j)=\rho(k) \right).
\end{equation}
Additionally, we denote $r_\mathrm{rate}:=\operatorname{rank}_\mathbb{Z}(\Lambda_\mathrm{rate})$ to be the $\mathbb{Z}$-rank of $\Lambda_\mathrm{rate}$. Whenever the scale assignment $\sigma$ is fixed, we suppress it from the notation and write $s$, $c_\mathrm{scale}$, $\Lambda_\mathrm{rate}$, and $r_\mathrm{rate}$. When comparing several admissible assignments, we instead write $s(\sigma)$, $c_\mathrm{scale}(\sigma)$, $\Lambda_\mathrm{rate}(\sigma)$, and $r_\mathrm{rate}(\sigma)$. When forming coordinate vectors, we relabel $U(\sigma)$ as $\{1,\dots,s(\sigma)\}$.
\subsection{The Valuation Homomorphism}\label{subsec:valuation}
Each active coefficient $c_{ij}$ is a Laurent monomial in the rate and concentration symbols, where negative exponents may occur due to the $1/K_i$ factor in \eqref{setcijeqn}. We therefore work in the Laurent polynomial ring
\begin{equation}
 R:= \mathbb{Z}\left[\kappa_1^{\pm 1}, \dots, \kappa_h^{\pm 1}, a_q^{\pm 1}\ (q\in U_\mathrm{cand})\right],
\end{equation}
and denote by
\begin{equation}
\mathcal{M}:=\left\{\prod_{p=1}^h\kappa_p^{\alpha_p}\prod_{q\in U_\mathrm{cand}} a_q^{\beta_q}\middle|\alpha_p,\beta_q\in\mathbb{Z}\right\}\leq R^\times
\end{equation}
the multiplicative group of sign-free Laurent monomials, and by $\mathcal{C}$ the multiplicative subgroup of $\mathcal{M}$ generated by all $c_{ij}$ with $(i,j)\in\mathcal{P}$. We define the valuation homomorphism
\begin{equation}
\operatorname{val}: \mathcal{M} \longrightarrow \mathbb{Z}^{h}\oplus \mathbb{Z}^{U_\mathrm{cand}},\qquad
\operatorname{val}\left(\prod_{p=1}^h \kappa_p^{\alpha_p}\prod_{q} a_q^{\beta_q}\right)=(\alpha;\beta), \label{valhomomorphism}
\end{equation}
which is an isomorphism of abelian groups (on all of $R^\times$ the corresponding map would have kernel $\{\pm1\}$).\\~\\
Let $f_p$ denote the $p$-th standard basis vector of $\mathbb{Z}^h$, let $e_q$ ($q\in U_\mathrm{cand}$) denote the standard basis vectors of $\mathbb{Z}^{U_\mathrm{cand}}$, and let $b_i:=e_{\sigma(i)}$ be the one-hot vector encoding the scale of species $i$, so that if $K_i=a_p$, then $b_i = e_p$. Then, for $(i,j)\in\mathcal{P}$, substituting \eqref{setcijeqn} into \eqref{valhomomorphism} leads to
\begin{equation}\label{eq:valcij}
    \operatorname{val}(c_{ij}) = \left(f_{\rho(j)};\ \sum_{\ell=1}^n \nu_{\ell j} b_\ell - b_i\right).
    \end{equation}
We define the \textbf{exponent matrix} $A$ to be the matrix with columns $\operatorname{val}(c_{ij})$ for $(i,j)\in\mathcal{P}$, and let $\Gamma = \operatorname{Im}_\mathbb{Z}(A)$ be the corresponding \textbf{exponent lattice}. Since only used symbols occur, $\Gamma\subseteq\mathbb{Z}^h\oplus\mathbb{Z}^{U}$, and we identify $\mathbb{Z}^{U}$ with $\mathbb{Z}^s$ where convenient. Under this identification, $\operatorname{val}(\mathcal{C})=\Gamma$, whose rank equals $\operatorname{rank}_{\mathbb{Z}}(A)$ and the number of nonzero invariant factors in the Smith normal form (SNF) of $A$.\\~\\
Finally, let $\deg_\kappa:\mathbb{Z}^h\oplus\mathbb{Z}^{U_\mathrm{cand}}\to\mathbb{Z}$ denote the sum of the rate coordinates (the \textbf{rate degree}), so that $\deg_\kappa(\operatorname{val}(c_{ij}))=1$ for all $(i,j)\in\mathcal{P}$, and let $\mathcal{D}:\mathbb{Z}^h\oplus\mathbb{Z}^{U_\mathrm{cand}}\to\mathbb{Z}^k$ be the \textbf{dimension homomorphism} with $\mathcal{D}(f_p;0)=\dim(\kappa_p)$ and $\mathcal{D}(0;e_q)=\dim(a_q)$, so that $\mathcal{D}(\operatorname{val}(g))$ is the dimension of $g\in\mathcal{M}$.
\begin{lemma}\label{lem:timedeg}
Under dimensional consistency, the time component of $\mathcal{D}(v)$ equals $-\deg_\kappa(v)$ for every $v\in\mathbb{Z}^h\oplus\mathbb{Z}^{U_\mathrm{cand}}$. In particular, $\deg_\kappa(v_T)=-1$ for every time scale $K_T$, where $v_T:=\operatorname{val}(K_T)$.
\end{lemma}
\begin{proof}
Every rate symbol occurs in a reaction with an active species, so \eqref{eq:dimconsistency} holds for some $(i,j)$ with $\rho(j)=p$, and since $\dim(X_\ell)$ has zero time component, the time component of $\dim(\kappa_p)$ is $-1$. Noting that concentration symbols have zero time component, both sides of the claimed identity are therefore homomorphisms that agree on the generators $(f_p;0)$ and $(0;e_q)$.
\end{proof}
In this work, we aim to determine, for a given non-dimensionalisation, the minimal number of multiplicative dimensionless parameters needed to express the coefficients $\{\tilde{c}_{ij}\}$ and the nonzero initial conditions \eqref{eq:dimlessIC} as Laurent monomials, and characterise which choices of species scales minimise this number.
\section{Results\label{resultssection}}
Before we state and prove the main theorems of this work, it is prudent to present two lattice-theoretic facts.
\begin{lemma}\label{graphranklemma}
Let $H$ be a finite undirected graph with vertex set $W$, edge set $E$ and connected components $D_1,\ldots,D_{c(H)}$. Then
\begin{equation}
\mathcal{L}(H):=\operatorname{span}_{\mathbb{Z}}\left(e_u-e_v\mid\{u,v\}\in E\right)
=\left\{\beta\in\mathbb{Z}^W\ \middle|\ \sum_{u\in D_\eta}\beta_u=0\text{ for every }\eta\right\}.
\end{equation}
Consequently, $\mathcal{L}(H)$ is saturated and has rank $|W|-c(H)$.
\end{lemma}
\begin{lemma}\label{rankadditivitylemma}
If $\Lambda\leq\mathbb{Z}^N$ and $\varphi:\mathbb{Z}^N\to\mathbb{Z}^d$ is a homomorphism, then
\begin{equation}
\operatorname{rank}_{\mathbb{Z}}(\Lambda)
=\operatorname{rank}_{\mathbb{Z}}(\varphi(\Lambda))
+\operatorname{rank}_{\mathbb{Z}}(\Lambda\cap\ker(\varphi)).
\end{equation}
\end{lemma}
These are standard results in algebraic graph theory and the theory of free abelian groups \cite{Godsil2001,Biggs1974,Lang2002}; however, for completeness, we provide proofs in Section B.
\begin{numberedremark}[Monomial relations among rate constants]\label{rem:rate-relations}
The independence assumption excludes, for example, the monomial constraints associated with detailed balance in reversible networks \cite{Wegscheider1901}. If homogeneous relations $\kappa^w=1$ are imposed for exponent vectors generating a lattice $\mathcal{W}\subseteq\mathbb{Z}^h\oplus\{0\}$, positivity implies the same relations for its saturation $\mathcal{W}_{\mathrm{sat}}=(\mathcal{W}\otimes_\mathbb{Z}\mathbb{Q})\cap\mathbb{Z}^{h+|U_\mathrm{cand}|}$. When these are the only algebraic relations imposed, monomial parameter lattices are replaced by their images in the quotient by $\mathcal{W}_{\mathrm{sat}}$. The rank of the image of a lattice $\Lambda$ is
\begin{equation}
\operatorname{rank}_\mathbb{Z}(\Lambda)-\operatorname{rank}_\mathbb{Z}(\Lambda\cap\mathcal{W}_{\mathrm{sat}}),
\end{equation}
by \Cref{rankadditivitylemma}. The independent-symbol counts are therefore upper bounds after such relations are imposed, with general algebraic relations requiring additional treatment.
\end{numberedremark}
\begin{theorem}\label{maintheorem}
Fix a scale assignment $\sigma$ for a CRN of the form \eqref{maincrnform}, with $h$ distinct rate symbols, $s$ used scale symbols and $c_\mathrm{scale}$ connected components of $H_\mathrm{scale}$. If $r_\mathrm{rate}$ is the rank of $\Lambda_\mathrm{rate}$, then the exponent lattice, $\Gamma$, generated by the columns of $A$ has rank
\begin{equation}
\operatorname{rank}_{\mathbb{Z}}(\Gamma)=h+s-c_\mathrm{scale}+r_\mathrm{rate}. \label{maineqn}
\end{equation}
\end{theorem}
\begin{proof}
Denote $\pi_\kappa: \mathbb{Z}^{h} \oplus \mathbb{Z}^{s} \to \mathbb{Z}^h$ as the projection operator onto rate-symbol coordinates, and $\pi_a: \mathbb{Z}^{h} \oplus \mathbb{Z}^{s}\to\mathbb{Z}^s$ as the projection operator onto scale-symbol coordinates. By \eqref{eq:valcij}, $\pi_\kappa(\operatorname{val}(c_{ij}))=f_{\rho(j)}$ for every active species-reaction pair $(i,j)$. Because each reaction has at least one active species, we have that
\begin{equation}
  \pi_\kappa(\Gamma) = \operatorname{span}_{\mathbb{Z}}(f_{\rho(j)} | j \in \set{1,\dots,r}) = \operatorname{span}_\mathbb{Z}(f_p | p\in \set{1,\dots,h}) \cong \mathbb{Z}^h,
\end{equation}
and so $\operatorname{rank}_{\mathbb{Z}}(\pi_\kappa(\Gamma)) = h$.\\~\\
Next, we define 
\begin{equation}
\mathcal{H}:=\Gamma\cap\ker(\pi_\kappa)\subseteq\{0\}^h\oplus\mathbb{Z}^{s}\cong\mathbb{Z}^s.
\end{equation}
By \Cref{rankadditivitylemma}, applied to $\Lambda=\Gamma$ and $\varphi = \pi_\kappa$, we have that
\begin{equation}
  \operatorname{rank}_\mathbb{Z}(\Gamma) = \operatorname{rank}_\mathbb{Z}(\pi_\kappa(\Gamma)) + \operatorname{rank}_\mathbb{Z}(\mathcal{H}) = h + \operatorname{rank}_\mathbb{Z}(\mathcal{H}). \label{rankzgammaeqn}
\end{equation}
For any reaction $j$ and any two active species $i,i^\prime \in\{1,\ldots,n\}$ with $P_{ij}=P_{i^\prime j}=1$, substituting in \eqref{eq:valcij} and subtracting leads to 
\begin{equation}
  \operatorname{val}(c_{i^\prime j}) - \operatorname{val}(c_{ij}) = \left(0; b_i-b_{i^\prime}\right) \in \mathcal{H}.
\end{equation}
We therefore define the scale-difference lattice
\begin{equation}
\mathcal{L}_{\mathrm{scale}} := \operatorname{span}_\mathbb{Z}\left(b_{i^\prime}-b_i \mid \set{i,i^\prime} \in E(G)\right)=\operatorname{span}_\mathbb{Z}\left(b_i-b_{i^\prime} \mid \set{i,i^\prime} \in E(G)\right) \subseteq \mathbb{Z}^s, \label{Lscaleformula}
\end{equation}
with the coefficient differences above giving $\{0\}^h\oplus\mathcal{L}_\mathrm{scale}\subseteq\mathcal{H}$. We can decompose \eqref{Lscaleformula} by components, so that
\begin{equation}
    \mathcal{L}_{\mathrm{scale}} = \sum_{\alpha = 1}^q \mathcal{L}_\alpha, \quad \text{where }\mathcal{L}_\alpha:= \operatorname{span}_\mathbb{Z}\left(b_{i^\prime}-b_i \mid \set{i,i^\prime} \in E(G), i, i^\prime \in C_\alpha \right).
\end{equation}
For any $\alpha\in \{1,\dots,q\}$, since $C_\alpha$ is connected, for any two species $u,v \in C_\alpha$, we have that $b_v - b_u \in \mathcal{L}_\alpha$. Since $b_i = e_{\sigma(i)}$, it follows that for any $m,p \in U_\alpha$, there exist $u,v \in C_\alpha$ with $\sigma(u)=m$ and $\sigma(v)=p$, so that $e_p - e_m = b_v - b_u \in \mathcal{L}_\alpha$. Conversely, $b_{i^\prime} - b_{i} = e_{\sigma(i^\prime)} - e_{\sigma(i)}$ with $\sigma(i^\prime), \sigma(i) \in U_\alpha$. Therefore,
\begin{equation}
\mathcal{L}_\alpha = \operatorname{span}_\mathbb{Z}\left(e_m-e_p \middle | m,p \in U_\alpha\right) = \mathcal{L}\left(K_{U_\alpha} \right),
\end{equation}
which is the difference lattice of the clique $K_{U_\alpha}$ on $U_\alpha$. Furthermore, we have that $\mathcal{L}_{\mathrm{scale}} = \sum_{\alpha=1}^q \mathcal{L}_\alpha$ is exactly the lattice $\mathcal{L}(H_\mathrm{scale})$, since $H_\mathrm{scale}$, by definition and \eqref{hscaledef}, is the union of cliques on $U_\alpha$, so that $E(H_{\mathrm{scale}})=\bigcup_{\alpha=1}^q E(K_{U_\alpha})$ and each $\mathcal{L}_\alpha=\mathcal{L}(K_{U_\alpha})$. Hence, by \Cref{graphranklemma}, we have that
\begin{equation}
 \operatorname{rank}_\mathbb{Z}(\mathcal{L}_{\mathrm{scale}}) = \operatorname{rank}_\mathbb{Z}(\mathcal{L}(H_\mathrm{scale})) = |U| - c_\mathrm{scale} = s - c_\mathrm{scale}. \label{lscalerankform}
\end{equation}
To capture the additional directions in $\mathcal{H}$ arising from reactions that share a rate constant but have different stoichiometries, we pass to the quotient group
\begin{equation}
    Q := \mathbb{Z}^s / \mathcal{L}_{\mathrm{scale}}.
\end{equation}
Let $D_1,\dots,D_{c_\mathrm{scale}}$ denote the connected components of $H_\mathrm{scale}$, and let $w_1,\dots,w_{c_\mathrm{scale}}$ be the standard basis of $\mathbb{Z}^{c_\mathrm{scale}}$. Define
\begin{equation}
\chi:\mathbb{Z}^s\to\mathbb{Z}^{c_\mathrm{scale}},\qquad \chi(e_p)=w_\eta \quad\text{if $p\in D_\eta$}.
\end{equation}
Since every generator of $\mathcal{L}_{\mathrm{scale}}=\mathcal{L}(H_\mathrm{scale})$ lies in $\ker(\chi)$, $\mathcal{L}_{\mathrm{scale}}\subseteq\ker(\chi)$. Conversely, letting $\beta=\sum_{p=1}^s\beta_pe_p\in\ker(\chi)$, we have that $\sum_{p\in D_\eta}\beta_p=0$ for each $\eta$. Choosing $p_\eta\in D_\eta$, we have
\begin{equation}
\beta=\sum_{\eta=1}^{c_\mathrm{scale}}\sum_{p\in D_\eta\setminus\set{p_\eta}}\beta_p(e_p-e_{p_\eta})\in\mathcal{L}_{\mathrm{scale}},
\end{equation}
since each $D_\eta$ is connected, so that $\ker(\chi)=\mathcal{L}_{\mathrm{scale}}$. Since $\chi$ is surjective, the first isomorphism theorem gives $Q\cong\mathbb{Z}^{c_\mathrm{scale}}$. Letting $\phi:\mathbb{Z}^s\to Q$ be the canonical quotient map, we identify $Q$ with this free abelian group and extend $\phi$ by ignoring the rate-symbol component to define
\begin{equation}
  \tilde{\phi}:\mathbb{Z}^h\oplus\mathbb{Z}^s\to Q,
  \qquad \tilde{\phi}(\alpha;\beta) := \phi(\beta).
\end{equation}
For each species $i\in\mathcal{I}$, this identification gives
\begin{equation}
  \phi(b_i) = w_{\psi(i)}.
\end{equation}
For each reaction $j\in\set{1,\dots,r}$, choose an active species $i_j$ with $P_{i_j j}=1$. Then
\begin{equation}
    z_j := \left(\sum_{\ell=1}^n \nu_{\ell j} w_{\psi(\ell)}\right) - w_{\psi_V(j)} = \phi\left(\sum_{\ell=1}^n \nu_{\ell j} b_\ell - b_{i_j}\right) \in Q.
\end{equation}
We again note that this is well-defined, since if $i_j^\prime$ is another active species with $P_{i_j^\prime j}=1$, then $b_{i_j^\prime} - b_{i_j}\in\mathcal{L}_\mathrm{scale}$, so that the corresponding expressions differ by an element of $\mathcal{L}_{\mathrm{scale}}$ and hence have the same image in $Q$. Furthermore, by the definition of $\Lambda_\mathrm{rate}$ and $z_j$, we have that $\Lambda_{\mathrm{rate}}\subset Q$.\\~\\
We now analyse the restriction $\tilde{\phi}|_{\mathcal{H}}$, and make the following claim:
\begin{claim}
We have that
\begin{equation}
   \operatorname{Im}_\mathbb{Z}(\tilde{\phi}|_{\mathcal{H}}) = \Lambda_{\mathrm{rate}}.
\end{equation}
\end{claim}
\begin{proof}
We prove this result by showing that $\Lambda_{\mathrm{rate}} \subseteq \operatorname{Im}_\mathbb{Z}(\tilde{\phi}|_{\mathcal{H}})$ and $\operatorname{Im}_\mathbb{Z}(\tilde{\phi}|_{\mathcal{H}}) \subseteq \Lambda_{\mathrm{rate}}$, so that $\operatorname{Im}_\mathbb{Z}(\tilde{\phi}|_{\mathcal{H}}) = \Lambda_{\mathrm{rate}}$.
\begin{enumerate}
    \item [(a)] If reactions $j$ and $k$ share the same reaction-rate symbol, so that $\kappa_{\rho(j)}=\kappa_{\rho(k)}$, then for active species $i_j$ in $j$ and $i_k$ in $k$,
\begin{equation}
    u_{jk} := \operatorname{val}(c_{i_j j}) - \operatorname{val}(c_{i_k k})\in\mathcal{H},
\end{equation}
and
\begin{equation}
    \tilde{\phi}(u_{jk}) = \phi\left(\sum_{\ell=1}^n \nu_{\ell j}b_\ell - b_{i_j}\right) - \phi\left(\sum_{\ell=1}^n \nu_{\ell k}b_\ell - b_{i_k}\right) = z_j - z_k.
\end{equation}
Thus, every generator $z_j - z_k$ of $\Lambda_{\mathrm{rate}}$ lies in $\operatorname{Im}_\mathbb{Z}(\tilde{\phi}|_{\mathcal{H}})$, so that
\begin{equation}
    \Lambda_{\mathrm{rate}} \subseteq \operatorname{Im}_\mathbb{Z}(\tilde{\phi}|_{\mathcal{H}}). \label{forwardinclusion}
\end{equation}
    \item [(b)] Conversely, any $u\in\mathcal{H}$ is an integer combination of valuation vectors, so that
\begin{equation}
    u = \sum_{\substack{1\leq i\leq n,\ 1\leq j\leq r\\ P_{ij}=1}} \lambda_{ij} \operatorname{val}(c_{ij}), \qquad \lambda_{ij}\in\mathbb{Z}.
\end{equation}
Denoting $\mu_j :=\sum\limits_{i \mid P_{ij}=1} \lambda_{ij}$, we have that since $\pi_\kappa(\operatorname{val}(c_{ij})) = f_{\rho(j)}$ for active pairs,
\begin{equation}
    \pi_\kappa(u) = \sum_{\substack{1\leq i\leq n,\ 1\leq j\leq r\\ P_{ij}=1}}\lambda_{ij} f_{\rho(j)} = \sum_{j=1}^r  \mu_j f_{\rho(j)}.
\end{equation}
The condition $u\in\ker(\pi_\kappa)$ is therefore equivalent to
\begin{equation}
    \sum_{j| \rho(j)=p} \mu_j = 0, \quad\text{for every $p\in\{1,\dots,h\}$}. \label{kernelcond}
\end{equation}
Moreover, for any active species $i_j$ in reaction $j$, we have that
\begin{equation}
  \operatorname{val}(c_{ij}) - \operatorname{val}(c_{i_j j}) = (0; b_{i_j} - b_i) \in \{0\}^h\oplus \mathcal{L}_{\mathrm{scale}},
\end{equation}
so $\tilde{\phi}(\operatorname{val}(c_{ij})) = \tilde{\phi}(\operatorname{val}(c_{i_j j})) = z_j$. Hence,
\begin{equation}
    \tilde{\phi}(u) = \sum_{\substack{1\leq i\leq n,\ 1\leq j\leq r\\ P_{ij}=1}} \lambda_{ij} \tilde{\phi}(\operatorname{val}(c_{ij})) = \sum_{j=1}^r \mu_j z_j.
\end{equation}
For each rate symbol $p\in\set{1,\dots,h}$, we choose an anchor reaction $j_p$ with $\rho(j_p)=p$, so that
\begin{equation}
    \sum_{j|\rho(j)=p} \mu_j z_j = \sum_{j | \rho(j)=p} \mu_j (z_j - z_{j_p}) + \left(\sum_{j| \rho(j)=p} \mu_j\right) z_{j_p} = \sum_{j| \rho(j)=p} \mu_j (z_j - z_{j_p}),
\end{equation}
since $\sum_{j|\rho(j)=p} \mu_j = 0$ by \eqref{kernelcond}. Therefore,
\begin{equation}
    \tilde{\phi}(u) = \sum_{p=1}^h \sum_{j|\rho(j)=p} \mu_j (z_j - z_{j_p}) \in \Lambda_{\mathrm{rate}},
\end{equation}
because each $z_j - z_{j_p}$ is one of the generators of $\Lambda_{\mathrm{rate}}$ with $\rho(j)=\rho(j_p)$.
Thus,
\begin{equation}
    \operatorname{Im}(\tilde{\phi}|_{\mathcal{H}}) \subseteq \Lambda_{\mathrm{rate}}. \label{reverseinclusion}
\end{equation}
\end{enumerate}
Combining \eqref{forwardinclusion} and \eqref{reverseinclusion} gives the required result.
\end{proof}
Additionally, by the definition of $Q$, we have that
\begin{equation}
    \ker_\mathbb{Z}(\tilde{\phi}|_{\mathcal{H}})
    = \{0\}^h\oplus\mathcal{L}_{\mathrm{scale}}.
\end{equation}
Thus, by \Cref{rankadditivitylemma}, we have that
\begin{equation}
    \operatorname{rank}_\mathbb{Z}(\mathcal{H}) = \operatorname{rank}_\mathbb{Z}(\mathcal{L}_{\mathrm{scale}}) + \operatorname{rank}_\mathbb{Z}(\Lambda_{\mathrm{rate}}) = s - c_\mathrm{scale} + r_{\mathrm{rate}}. \label{rankzhequation}
\end{equation}
Thus, substituting \eqref{rankzhequation} into \eqref{rankzgammaeqn} results in
\begin{equation}
    \operatorname{rank}_\mathbb{Z}(\Gamma) = h + s - c_\mathrm{scale} + r_{\mathrm{rate}},   
\end{equation}
as desired.
\end{proof}
It is important to note that the quantity $r_\mathrm{rate}$ is attached to a fixed scale assignment and can change when the assignment changes, with this dependence characterised in \Cref{rassignedprop} through a lattice determined by the reaction data alone.
\begin{proposition}\label{rassignedprop}\label{proof:rassignedprop}
Let $\tau:\mathcal{I}\to\{1,\ldots,q\}$ map each species to its component of $G$, let $\tau_V(j)$ be the common component of the active species in reaction $j$, and let $u_1,\ldots,u_q$ be the standard basis of $\mathbb{Z}^q$. Define
\begin{equation}
\bar{z}_j:=\sum_{\ell=1}^n\nu_{\ell j}u_{\tau(\ell)}-u_{\tau_V(j)},
\qquad
\bar{\Lambda}_\mathrm{rate}:=\operatorname{span}_{\mathbb{Z}}(\bar{z}_j-\bar{z}_k\mid\rho(j)=\rho(k)).
\end{equation}
For an assignment $\sigma$, let $\eta_\sigma(\alpha)$ be the component of $H_\mathrm{scale}$ containing $U_\alpha$, and define $\pi_\sigma(u_\alpha)=w_{\eta_\sigma(\alpha)}$. Then $z_j(\sigma)=\pi_\sigma(\bar{z}_j)$, $\Lambda_\mathrm{rate}(\sigma)=\pi_\sigma(\bar{\Lambda}_\mathrm{rate})$, and
\begin{equation}
r_\mathrm{rate}(\sigma)=\operatorname{rank}_{\mathbb{Z}}(\bar{\Lambda}_\mathrm{rate})- \operatorname{rank}_{\mathbb{Z}}(\bar{\Lambda}_\mathrm{rate}\cap\ker(\pi_\sigma)). \label{rrateassignmentformula}
\end{equation}
Consequently, $r_\mathrm{rate}(\sigma)$ is constant on an admissible family $\Sigma$ if and only if $\operatorname{rank}_\mathbb{Z}\left(\bar\Lambda_\mathrm{rate}\cap\ker(\pi_\sigma)\right)$ is constant on $\Sigma$. In particular, this holds if $\ker(\pi_\sigma)$ is fixed or if $\bar{\Lambda}_\mathrm{rate}=\{0\}$.
\end{proposition}
\begin{proof}
For each $\alpha$, the set $U_\alpha$ is a clique of $H_\mathrm{scale}$, so it lies in a single connected component and $\eta_\sigma(\alpha)$ is well-defined. Since each species $\ell$ in $C_\alpha$ has $\sigma(\ell)\in U_\alpha$, we have $\psi(\ell)=\eta_\sigma(\tau(\ell))$; likewise, if $i$ is any active species in reaction $j$, then $\psi_V(j)=\eta_\sigma(\tau_V(j))$. Substituting these identities into the definition of $z_j$ gives $z_j(\sigma)=\pi_\sigma(\bar{z}_j)$, and taking within-rate differences leads to $\Lambda_\mathrm{rate}(\sigma)=\pi_\sigma(\bar{\Lambda}_\mathrm{rate})$. Thus, $r_\mathrm{rate}(\sigma)=\operatorname{rank}_\mathbb{Z}\left(\pi_\sigma(\bar{\Lambda}_\mathrm{rate})\right)$.\\~\\
Applying \Cref{rankadditivitylemma} to the subgroup $\bar{\Lambda}_\mathrm{rate}\leq\mathbb{Z}^q$ and the homomorphism $\pi_\sigma$ yields \eqref{rrateassignmentformula}. Since $\operatorname{rank}_\mathbb{Z}(\bar{\Lambda}_\mathrm{rate})$ is assignment-independent, we have that $r_\mathrm{rate}(\sigma)$ is independent of $\sigma\in\Sigma$ if and only if $\operatorname{rank}_\mathbb{Z}\left(\bar{\Lambda}_\mathrm{rate}\cap\ker(\pi_\sigma)\right)$ is independent of $\sigma\in\Sigma$. We note that the two stated conditions are sufficient: if $\ker(\pi_\sigma)$ is fixed on $\Sigma$, or if the lattice being projected is already zero, then $\operatorname{rank}_\mathbb{Z}\left(\bar{\Lambda}_\mathrm{rate}\cap\ker(\pi_\sigma)\right)$ is assignment-independent.
\end{proof}
\begin{corollary}\label{rratecorollary}
If all reaction-rate symbols are distinct, then $\bar{\Lambda}_\mathrm{rate}=\{0\}$ and $\operatorname{rank}_{\mathbb{Z}}(\Gamma)=h+s-c_\mathrm{scale}$. If $G$ is connected, then $c_\mathrm{scale}=1$ and $z_j=(\sum_\ell\nu_{\ell j}-1)w_1$. In this case $r_\mathrm{rate}\in\{0,1\}$ is assignment-independent, with value one precisely when two reactions sharing a rate symbol have different total reactant orders. Thus, $\operatorname{rank}_{\mathbb{Z}}(A)=h+s+r_\mathrm{rate}-1$ depends on $\sigma$ only through $s$, and is unchanged by rearrangements preserving the multiset of scale symbols.
\end{corollary}
\begin{proof}
Clearly, distinct rate symbols give no nonzero within-rate differences so that $\bar{\Lambda}_\mathrm{rate}=\{0\}$. If $G$ is connected, the scale graph $H_{\mathrm{scale}}$ is the clique on $U$ and $\pi_\sigma:\mathbb{Z}\to\mathbb{Z}$ is the identity, giving the stated formula for $z_j$ and its consequences.
\end{proof}
It follows that to minimise $\operatorname{rank}_\mathbb{Z}(\Gamma)$, and consequently $\operatorname{rank}_\mathbb{Z}(\mathcal{C})$, we must minimise $s-c_\mathrm{scale}+r_\mathrm{rate}$ over admissible scale assignments. On any admissible family $\Sigma$ for which \Cref{rassignedprop} shows that $r_\mathrm{rate}(\sigma)$ is constant, this reduces to minimising $s-c_\mathrm{scale}$ over $\sigma\in\Sigma$.

\begin{numberedremark}[Dimensional interpretation of $r_\mathrm{rate}$]\label{rem:rrate}
If all species have a common nonzero dimension, \eqref{eq:dimconsistency} forces reactions sharing a rate symbol to have equal total reactant order, so that every generator of $\Lambda_\mathrm{rate}$ has coordinate sum zero, giving $r_\mathrm{rate}=0$ for connected $G$. For disconnected networks this rank may nevertheless be positive, as the mutual-activation system in \Cref{ex:rrate} demonstrates. The nonzero qualification is essential, since dimensionless species may occur with different orders in reactions sharing the same rate symbol.
\end{numberedremark}

\begin{example}[A network with $r_\mathrm{rate}=1$]\label{ex:rrate}
Consider the reactions
\begin{equation}
\mathrm{X}_2\xrightarrow{\kappa}\mathrm{X}_2+\mathrm{X}_1,
\qquad
\mathrm{X}_1\xrightarrow{\kappa}\mathrm{X}_1+\mathrm{X}_2,
\qquad
\mathrm{X}_1\xrightarrow{\kappa_1}\emptyset,
\qquad
\mathrm{X}_2\xrightarrow{\kappa_2}\emptyset,
\end{equation}
where both $X_1$ and $X_2$ have the same dimension. The governing equations are $dX_1/dT=\kappa X_2-\kappa_1X_1$ and $dX_2/dT=\kappa X_1-\kappa_2X_2$, with $h=3$, and each reaction has a single active species, so $G$ consists of two isolated vertices. For distinct scale symbols $K_1=a$ and $K_2=b$, the scale graph also has two isolated vertices, while $z_1=w_2-w_1$ and $z_2=w_1-w_2$, so that $\Lambda_\mathrm{rate}=2\mathbb{Z}(w_1-w_2)$, $r_\mathrm{rate}=1$, and \Cref{maintheorem} gives $\operatorname{rank}_\mathbb{Z}(\Gamma)=4$. With $K_T=\kappa_1^{-1}$, the dimensionless coefficients are
\begin{equation}
1,\qquad \frac{\kappa_2}{\kappa_1},\qquad
\frac{\kappa b}{\kappa_1a},\qquad
\frac{\kappa a}{\kappa_1b},
\end{equation}
and their corresponding lattice has rank three, with the ratio and product of the two activation coefficients being $(b/a)^2$ and $(\kappa/\kappa_1)^2$, respectively, which are algebraically independent under the standing assumptions. If we were to choose $K_1=K_2=a$ instead, then $c_\mathrm{scale}=1$, $r_\mathrm{rate}=0$ and $\operatorname{rank}_\mathbb{Z}(\Gamma)=3$.
\end{example}

\subsection{Existence and Unimodular Uniqueness of Laurent Bases}
We now formalise what it means to generate all coefficients from a set of monomials and then prove unimodular uniqueness of the bases in the polynomial and Laurent polynomial settings.
\begin{definition}\label{Laurentdefs}
Let $\mathcal{C}_0\subset\mathcal{M}$ be finite and let $\Lambda$ be the lattice generated by its valuations. A finite set $\mathcal{G}=\{g_1,\ldots,g_m\}\subset\mathcal{M}$ is a \textbf{Laurent generating set} for $\mathcal{C}_0$ if every $c\in\mathcal{C}_0$ has an expression $c=\prod_{\ell=1}^m g_\ell^{\gamma_\ell}$ with $\gamma_\ell\in\mathbb{Z}$, and a \textbf{polynomial generating set} if these exponents can be taken in $\mathbb{N}$. It is a \textbf{Laurent basis}, respectively a \textbf{polynomial basis}, if its valuations additionally form a $\mathbb{Z}$-basis of $\Lambda$.
\end{definition}
We apply these definitions to the coefficient lattices $\Gamma$ and $\tilde{\Gamma}$, and subsequently to the coefficients together with the nonzero initial conditions.
\begin{proposition}\label{laurentbasisexistence}
For any finite $\mathcal{C}_0\subset\mathcal{M}$ with exponent lattice $\Lambda$,
\begin{enumerate}
\item[(a)] Laurent bases exist and have minimum cardinality $\operatorname{rank}_{\mathbb{Z}}(\Lambda)$ among Laurent generating sets;
\item[(b)] if every element of $\mathcal{C}_0$ is dimensionless, every element of every Laurent or polynomial basis is dimensionless.
\end{enumerate}
\end{proposition}
\begin{proof}
Clearly, the subgroup $\Lambda$ of a free abelian group is free, and the inverse images of any of its integral bases under $\operatorname{val}$ give a Laurent basis, and since every generating set spans a lattice containing $\Lambda$, its cardinality is at least $\operatorname{rank}_{\mathbb{Z}}(\Lambda)$. Furthermore, if $\mathcal{C}_0$ is dimensionless, then $\Lambda\subseteq\ker(\mathcal{D})$, which also contains every basis valuation, as required.
\end{proof}
\begin{theorem}\label{unimodularmain}
Let $\mathcal{G}$ and $\mathcal{G}^\prime$ be two bases, each Laurent or polynomial, of the same finite monomial set, with valuation matrices $W$ and $W^\prime$. If the generated lattice has rank $m$, then there is a unique $U\in\operatorname{GL}_m(\mathbb{Z})$ such that
\begin{equation}
W^\prime=WU,\qquad g^\prime_\ell=\prod_{\lambda=1}^m g_\lambda^{U_{\lambda\ell}}. \label{unishift}
\end{equation}
\end{theorem}
\begin{proof}
As $\mathcal{G}$ and $\mathcal{G}^\prime$ define isomorphisms from $\mathbb{Z}^m$ to the same lattice, composing one with the inverse of the other gives the required unimodular automorphism, and \eqref{unishift} follows from the injectivity of $\operatorname{val}$.
\end{proof}
We note that minimum cardinality alone does not imply this uniqueness as, for example, $\{\zeta\}$ and $\{\zeta^2\}$ both generate $\{\zeta^2\}$ when $\operatorname{val}(\zeta)\neq 0$, but only the latter is a basis of its exponent lattice.

\subsubsection{Time Scales and the Dimensionless Coefficient Lattice}\label{subsec:timescale}
Recalling that $\tilde{c}_{ij}=K_Tc_{ij}$ and $v_T=\operatorname{val}(K_T)$, we write $\tilde{\Gamma}=\operatorname{span}_{\mathbb{Z}}(\operatorname{val}(\tilde{c}_{ij})\mid(i,j)\in\mathcal{P})$, $\Gamma_\mathbb{Q}=\operatorname{span}_{\mathbb{Q}}(\Gamma)$ and $d=\operatorname{rank}_{\mathbb{Z}}(\Gamma)$, and define
\begin{equation}
\Delta\Gamma:=\Gamma\cap\ker(\deg_\kappa). \label{deltagammaeqn}
\end{equation}
Every coefficient valuation has rate degree one, so $\deg_\kappa(\Gamma)=\mathbb{Z}$ and $\operatorname{rank}_{\mathbb{Z}}(\Delta\Gamma)=d-1$.
\begin{lemma}\label{Gammalatticelemma}
For $f(v):=v+\deg_\kappa(v)v_T$ on the ambient exponent lattice,
\begin{enumerate}
\item[(a)] $\tilde{\Gamma}=f(\Gamma)$;
\item[(b)] for an arbitrary monomial valuation $v_T$, the rank of $\tilde{\Gamma}$ is $d-1$ if $v_T\in\Gamma_\mathbb{Q}$ and $\deg_\kappa(v_T)=-1$, and is $d$ otherwise;
\item[(c)] for a dimensionally admissible time scale, $\deg_\kappa(v_T)=-1$, so $\operatorname{rank}_{\mathbb{Z}}(\tilde{\Gamma}) = d-1$ if and only if $v_T\in\Gamma_\mathbb{Q}$;
\item[(d)] if $v_T\in\Gamma$ and $K_T$ has the dimension of time, then $\tilde{\Gamma}=\Delta\Gamma$. For an \textbf{anchor} $K_T=c_{i_0j_0}^{-1}$, this lattice is generated by the differences $\operatorname{val}(c_{ij})-\operatorname{val}(c_{i_0j_0})$.
\end{enumerate}
\end{lemma}
\begin{proof}
Since $\deg_\kappa$ is a homomorphism, so is $f$. For every active pair $(i,j)$, the coefficient valuation has rate degree one, and hence
  \begin{equation*}
  f(\operatorname{val}(c_{ij}))=\operatorname{val}(c_{ij})+v_T=\operatorname{val}(\tilde{c}_{ij}).
  \end{equation*}
  These coefficient valuations generate $\Gamma$, so their images generate $f(\Gamma)$, proving (a).
  
  For (b), letting $\gamma\in\Gamma$ and setting $m=\deg_\kappa(\gamma)$, we have
  \begin{equation*}
  f(\gamma)=0\quad\Longleftrightarrow\quad\gamma=-m v_T,
  \end{equation*}
  so that if $f(\gamma)=0$ and $\gamma\neq 0$, then $m\neq0$, so $v_T=-\gamma/m\in\Gamma_\mathbb{Q}$. Taking rate degrees gives $m=-m\deg_\kappa(v_T)$, and therefore $\deg_\kappa(v_T)=-1$. Conversely, if these two conditions hold, clearing denominators gives a positive integer $q_0$ with $q_0v_T\in\Gamma$, with $v_T \neq 0$ as $\deg_\kappa(v_T)=-1$, and
  \begin{equation*}
  f(q_0v_T)=q_0\bigl(1+\deg_\kappa(v_T)\bigr)v_T=0.
  \end{equation*}
  Moreover, the kernel equation shows that every element of $\Gamma\cap\ker(f)$ is a multiple of $v_T$, so that this kernel has rank one precisely under the stated conditions and is trivial otherwise. By \Cref{rankadditivitylemma} and part (a),
  \begin{equation*}
  \operatorname{rank}_{\mathbb{Z}}(\tilde{\Gamma})
  =d-\operatorname{rank}_{\mathbb{Z}}(\Gamma\cap\ker(f)),
  \end{equation*}
  which proves (b). Part (c) immediately follows from \Cref{lem:timedeg}.
  
  For (d), admissibility gives $\deg_\kappa(v_T)=-1$. If $v_T\in\Gamma$ as well, then $f(\gamma)\in\Gamma$ and
  \begin{equation*}
  \deg_\kappa(f(\gamma))
  =\bigl(1+\deg_\kappa(v_T)\bigr)\deg_\kappa(\gamma)=0
  \end{equation*}
  for every $\gamma\in\Gamma$, so $f(\Gamma)\subseteq\Delta\Gamma$. Conversely, $f(\delta)=\delta$ for every $\delta\in\Delta\Gamma$, giving $\Delta\Gamma\subseteq f(\Gamma)$ and hence $\tilde{\Gamma}=f(\Gamma)=\Delta\Gamma$. For an anchor, $v_T=-\operatorname{val}(c_{i_0j_0})$, so part (a) gives the stated difference generators.
  \end{proof}
\begin{theorem}\label{laurentbasistheorem}
For an admissible assignment $\sigma$, let
\begin{equation}
N_\mathrm{coef}(\sigma):=h+s(\sigma)-c_\mathrm{scale}(\sigma)+r_\mathrm{rate}(\sigma)-1.
\end{equation}
The minimum cardinality of a Laurent generating set of $\{\tilde{c}_{ij}\}$ is $N_\mathrm{coef}(\sigma)$ when $v_T\in\Gamma_\mathbb{Q}$, and is $N_\mathrm{coef}(\sigma)+1$ otherwise. Every anchor
\begin{equation}
K_T=\frac{K_{i_0}}{\kappa_{\rho(j_0)}\prod_{\ell=1}^nK_\ell^{\nu_{\ell j_0}}}
=c_{i_0j_0}^{-1} \label{KTformula}
\end{equation}
attains the smaller count, and its Laurent bases are the monomial sets whose valuations form a $\mathbb{Z}$-basis of $\Delta\Gamma$.
\end{theorem}
\begin{proof}
Combine \Cref{maintheorem,Gammalatticelemma} with \Cref{laurentbasisexistence}(a).
\end{proof}
\subsection{Existence and Unimodular Uniqueness of Polynomial Bases}\label{subsec:polynomial-bases}
We now restrict to polynomial generators, where monomial coefficients must be polynomials in the generator coordinates.  For each component $C_\alpha$ of $G$, choose an anchor scale index $p_\alpha^\ast\in U_\alpha$ and define $R_{\alpha,p}:=a_{p_\alpha^\ast}/a_p$, with $R_{\alpha,p_\alpha^\ast}=1$. If $\psi_G(j)$ denotes the component containing the active species of reaction $j$, then the anchored reaction monomial and its coefficient factorisation are
\begin{align}
G_j^\ast&:=\kappa_{\rho(j)}\left(\prod_{\ell=1}^n a_{\sigma(\ell)}^{\nu_{\ell j}}\right)a_{p^\ast_{\psi_G(j)}}^{-1}, \notag\\
c_{ij}&=G_j^\ast R_{\psi_G(j),\sigma(i)}. \label{polyfactor}
\end{align}
\begin{proposition}\label{ganchtheorem}
The set $\mathcal{G}_\mathrm{anch}:=\{G_j^\ast\mid j=1,\ldots,r\}\cup\{R_{\alpha,p}\mid \alpha=1,\ldots,q,\ p\in U_\alpha\setminus\{p_\alpha^\ast\}\}$ is a polynomial generating set for $\{c_{ij}\}$ whose valuations span exactly $\Gamma$. Moreover, a polynomial basis exists with $d=h+s-c_\mathrm{scale}+r_\mathrm{rate}$ elements.
\end{proposition}
\begin{proof}
From \eqref{polyfactor}, it is immediate that $\mathcal{G}_{\mathrm{anch}}$ is a polynomial generating set for $\{c_{ij}\}$. Additionally, from \eqref{valhomomorphism}, we have that
\begin{align}
    \operatorname{val}(G_j^\ast) &= \left(f_{\rho(j)}; \sum_{\ell=1}^n \nu_{\ell j}b_\ell - e_{p_{\psi_G(j)}^\ast}\right), \\*
    \operatorname{val}(R_{\alpha,p})&=\left(0;e_{p_\alpha^\ast}-e_p\right).
\end{align}
Since every $c_{ij}$ is a product of elements of $\mathcal{G}_{\mathrm{anch}}$, the valuation vectors of $\mathcal{G}_{\mathrm{anch}}$ span a lattice containing $\Gamma$. For any edge $\set{i,i^\prime}\in E(G)$, and any reaction $j$ with $P_{ij}=P_{i^\prime j}=1$, we have $\operatorname{val}(c_{ij})-\operatorname{val}(c_{i^\prime j})=(0;b_{i^\prime}-b_i)\in\Gamma$, so $\{0\}^h\oplus\mathcal{L}_{\mathrm{scale}}\subseteq\Gamma$. Furthermore, as $p_\alpha^\ast,p\in U_\alpha$ and $K_{U_\alpha}$ is a clique in $H_\mathrm{scale}$, we have $e_{p_\alpha^\ast}-e_p\in\mathcal{L}_{\mathrm{scale}}$, and hence $\operatorname{val}(R_{\alpha,p})\in \Gamma$. Therefore, for any active species $i$ in reaction $j$, $\operatorname{val}(G_j^\ast)=\operatorname{val}(c_{ij})-\operatorname{val}(R_{\psi_G(j),\sigma(i)})\in\Gamma$. Thus, the valuation vectors of $\mathcal{G}_{\mathrm{anch}}$ span exactly $\Gamma$.\\~\\
To prove that a polynomial basis exists, let $\deg_\kappa:\Gamma\to\mathbb{Z}$ denote the homomorphism obtained by summing the rate-symbol coordinates, and choose an active pair $(i_0,j_0)$, setting $u_1:=\operatorname{val}(c_{i_0j_0})$. Since $\deg_\kappa(u_1)=1$, we have that $\Gamma=\mathbb{Z}u_1\oplus\ker(\deg_\kappa)$. Choosing a $\mathbb{Z}$-basis $u_2,\dots,u_d$ of $\ker(\deg_\kappa)$, every active coefficient has a unique expression
\begin{equation}
    \operatorname{val}(c_{ij})=u_1+\sum_{\ell=2}^d a_{\ell,ij}u_\ell, \qquad a_{\ell,ij}\in\mathbb{Z},
\end{equation}
since every active coefficient has rate-degree $1$. For each $\ell=2,\dots,d$, we set $m_\ell:=\min_{P_{ij}=1}a_{\ell,ij}$ and define $u_1^\prime:=u_1+\sum_{\ell=2}^d m_\ell u_\ell$, so that $\{u_1^\prime,u_2,\dots,u_d\}$ is also a $\mathbb{Z}$-basis of $\Gamma$. Letting $g_1,\dots,g_d\in\mathcal{M}$ be the unique monomials whose valuation vectors are $u_1^\prime,u_2,\dots,u_d$, the injectivity of $\operatorname{val}$ gives
\begin{equation}
    c_{ij}=g_1\prod_{\ell=2}^d g_\ell^{a_{\ell,ij}-m_\ell},
\end{equation}
where every exponent is non-negative, and hence $\{g_1,\dots,g_d\}$ is a polynomial basis for $\{c_{ij}\}$, with the formula for $d$ following from \Cref{maintheorem}.
\end{proof}
In practice, it is much easier to convert dimensionless Laurent coefficients into polynomial coefficients than to directly construct a polynomial generating set. We now show that polynomial generating sets can always be obtained from Laurent ones via a uniform coordinate-wise shift of exponent vectors when the common factor is dimensionless, and that, after anchor normalisation, the same shift also yields a polynomial basis from any Laurent basis.
\begin{lemma}\label{lem:global-shift}
Let $\mathcal{G}=\{g_1,\ldots,g_{N_g}\}$ be a Laurent generating set for $\{\tilde{c}_{ij}\}$, with coefficient exponent vectors $v^{(1)},\ldots,v^{(N)}\in\mathbb{Z}^{N_g}$. Set $M_i:=\min_kv_i^{(k)}$ and $M=(M_i)$. Then
\begin{equation}
v^{\prime(k)}:=v^{(k)}-M\in\mathbb{N}^{N_g},\qquad
g^{v^{(k)}}=g^M g^{v^{\prime(k)}},
\end{equation}
where $g^v:=\prod_i g_i^{v_i}$ and the factor $g^M$ is independent of $k$.
\end{lemma}
\begin{proof}
  We immediately see, by construction, that ${v^\prime}^{(k)}_i = v^{(k)}_i - M_i \geq 0$ for all $i = 1,\dots,N_g$ and $k=1,\dots,N$. By definition, we have that 
  \begin{equation}
      g^{v^{(k)}}
      = \prod_{i=1}^{N_g} g_i^{M_i + {v^\prime}^{(k)}_i}
      = \left(\prod_{i=1}^{N_g} g_i^{M_i}\right)
        \left(\prod_{i=1}^{N_g} g_i^{{v^\prime}^{(k)}_i}\right)
      = g^M g^{{v^\prime}^{(k)}},
  \end{equation}
  where $g^M:= \prod_{i=1}^{N_g} g_i^{M_i}$ is independent of $k$. Thus, if $g^M$ is dimensionless, then, since this factor appears uniformly in every coefficient, we may absorb it into the time scale by replacing $K_T\mapsto K_T g^{-M}\in \mathcal{M}$, and in the resulting non-dimensional system all $\tilde{c}_{ij}$ are polynomials in the $g_\ell$.
\end{proof}
\begin{theorem}\label{thm:poly-shift}
Under the hypotheses of \Cref{lem:global-shift}, suppose that $g^M$ is dimensionless. After the time-scale transformation
\begin{equation}
K_T\mapsto K_Tg^{-M}, \label{KTformula2}
\end{equation}
$\mathcal{G}$ is a polynomial generating set for the rescaled coefficients $\{\tilde{c}_{ij}\}$. Furthermore, if $\mathcal{G}$ is a Laurent basis and the coefficient exponent vectors contain zero, then this procedure yields a polynomial basis for the rescaled coefficient lattice, and the difference lattice $\Delta\Gamma$ is unchanged. In particular, an anchor-normalised Laurent basis always satisfies these conditions, including dimensionlessness of $g^M$.
\end{theorem}
\begin{proof}
The transformation \eqref{KTformula2} multiplies every coefficient by $g^{-M}$, giving nonnegative exponents by \Cref{lem:global-shift} without changing the dimension of $K_T$. In particular, common translation leaves all valuation differences unchanged, and these differences generate $\Delta\Gamma$ by \Cref{Gammalatticelemma}. If zero is among the original exponent vectors, then $-M$ is among the shifted vectors and every original vector is recovered as $v^{(k)}=(v^{(k)}-M)-(-M)$. For a basis, the original vectors generate a lattice containing $\mathbb{Z}^{N_g}$, which, as they are integer vectors, is exactly $\mathbb{Z}^{N_g}$, so the shifted vectors generate the same lattice. An anchor supplies the zero vector, while \Cref{laurentbasisexistence}(b) establishes the required dimensionlessness.
\end{proof}
\begin{corollary}\label{laurentcorollary}
After anchor normalisation and, for polynomial bases, the transformation \eqref{KTformula2}, every Laurent or polynomial basis of the resulting coefficient set has $N_\mathrm{coef}(\sigma)$ elements. This is $h+s+r_\mathrm{rate}-2$ for connected $G$, and $h+s-c_\mathrm{scale}-1$ when all rate symbols are distinct, and any two such bases are related by a unimodular transformation.
\end{corollary}

\begin{theorem}\label{mingentheorem}
Let $\Sigma$ be a nonempty admissible family on which $r_\mathrm{rate}(\sigma)$ is constant, with common value $r_\mathrm{rate}$, and suppose that $\Sigma$ contains an assignment constant on each component of $G$. With the time scales of \Cref{laurentcorollary}, every coefficient basis associated with $\sigma\in\Sigma$ has at least $h+r_\mathrm{rate}-1$ elements, with this minimum being attained when $|U_\alpha(\sigma)|=1$ for every $\alpha$. For connected $G$, this is equivalent to all species scales coinciding.
\end{theorem}
\begin{proof}
The count is $h+s-c_\mathrm{scale}+r_\mathrm{rate}-1$, and $s\geq c_\mathrm{scale}$ with equality exactly when $H_\mathrm{scale}$ has no edges. Since this graph is the union of cliques on $U_\alpha$, equality holds exactly when every component of $G$ uses a single symbol, with the attainability assumption ensuring that such an assignment belongs to $\Sigma$.
\end{proof}

\subsection{Initial Conditions and the Parameter Count of the Full Initial Value Problem}\label{subsec:ivp}
For each component of $G$, define
\begin{equation}
W_\alpha:=\iota(\mathcal{J}\cap C_\alpha),\qquad
V_\alpha:=U_\alpha\cup W_\alpha,\qquad
U_\mathrm{tot}:=U\cup\iota(\mathcal{J}),
\end{equation}
and write $s_\mathrm{tot}=|U_\mathrm{tot}|$. Define the \textbf{augmented scale graph} as $H_\mathrm{aug}:=\bigcup_\alpha K_{V_\alpha}$ on $U_\mathrm{tot}$, with $c_\mathrm{aug}$ components, and let $\pi_\sigma^\mathrm{aug}:\mathbb{Z}^q\to\mathbb{Z}^{c_\mathrm{aug}}$ send $u_\alpha$ to the basis vector of the component containing $V_\alpha$, setting $r_\mathrm{IVP}:=\operatorname{rank}_{\mathbb{Z}}(\pi_\sigma^\mathrm{aug}(\bar{\Lambda}_\mathrm{rate}))$. The initial-condition and combined lattices are
\begin{equation}
\begin{aligned}
\mathcal{L}_\mathrm{IC}&:=\operatorname{span}_{\mathbb{Z}}\left((0;e_{\iota(i)}-e_{\sigma(i)})\mid i\in\mathcal{J}\right),\\
\Lambda_\mathrm{IVP}&:=\Gamma+\mathcal{L}_\mathrm{IC},\qquad
\tilde{\Gamma}_\mathrm{IVP}:=\tilde{\Gamma}+\mathcal{L}_\mathrm{IC},
\end{aligned}
\end{equation}
so that $\operatorname{rank}_{\mathbb{Z}}(\tilde{\Gamma}_\mathrm{IVP})$ is the minimum cardinality of a Laurent generating set for all dimensionless coefficients and nonzero initial conditions.

\begin{theorem}\label{thm:ivpcount}\label{proof:ivpcount}
For an admissible assignment $\sigma$, define
\begin{equation}
N_\mathrm{IVP}(\sigma):=h+s_\mathrm{tot}(\sigma)-c_\mathrm{aug}(\sigma)+r_\mathrm{IVP}(\sigma)-1. \label{eq:NIVP}
\end{equation}
For every dimensionally admissible time scale,
\begin{equation}
\operatorname{rank}_{\mathbb{Z}}(\tilde{\Gamma}_\mathrm{IVP})=
\begin{cases}
N_\mathrm{IVP}(\sigma),&v_T\in\Lambda_\mathrm{IVP}\otimes_{\mathbb{Z}}\mathbb{Q},\\
N_\mathrm{IVP}(\sigma)+1,&v_T\notin\Lambda_\mathrm{IVP}\otimes_{\mathbb{Z}}\mathbb{Q},
\end{cases} \label{eq:IVPtimecount}
\end{equation}
with every anchor attaining the smaller count. If $v_T\in\Lambda_\mathrm{IVP}$, then the integral lattices satisfy
\begin{equation}
\tilde{\Gamma}_\mathrm{IVP}=\Lambda_\mathrm{IVP}\cap\ker(\deg_\kappa)
=\Delta\Gamma+\mathcal{L}_\mathrm{IC}. \label{eq:IVPanchorlattice}
\end{equation}
Furthermore, if $W_\alpha\subseteq U_\alpha$ for all $\alpha$, in particular when $\sigma(i)=\iota(i)$ for every $i\in\mathcal{J}$, then $N_\mathrm{IVP}=N_\mathrm{coef}$.
\end{theorem}
\begin{proof}
Let $\pi_\kappa$ denote the rate-coordinate projection on $\mathbb{Z}^h\oplus\mathbb{Z}^{U_\mathrm{tot}}$. Since $\pi_\kappa(\Gamma)=\mathbb{Z}^h$ and $\mathcal{L}_\mathrm{IC}\subseteq\ker(\pi_\kappa)$, we have
\begin{equation}
\pi_\kappa(\Lambda_\mathrm{IVP})=\mathbb{Z}^h,
\qquad
\Lambda_\mathrm{IVP}\cap\ker(\pi_\kappa)
=\mathcal{H}+\mathcal{L}_\mathrm{IC},
\end{equation}
where $\mathcal{H}=\Gamma\cap\ker(\pi_\kappa)$ is as in the proof of \Cref{maintheorem}. It remains to compute the rank of the latter lattice, whose rate-coordinate block is zero.\\~\\
Letting $\chi_\mathrm{aug}:\{0\}^h\oplus\mathbb{Z}^{U_\mathrm{tot}}\to\mathbb{Z}^{c_\mathrm{aug}}$ send $(0;e_p)$ to the basis vector of the component of $H_\mathrm{aug}$ containing $p$, by \Cref{graphranklemma},
\begin{equation}
\ker(\chi_\mathrm{aug})=\{0\}^h\oplus\mathcal{L}(H_\mathrm{aug}),
\qquad
\operatorname{rank}_\mathbb{Z}(\mathcal{L}(H_\mathrm{aug}))
=s_\mathrm{tot}-c_\mathrm{aug}.
\end{equation}
The sublattice $\{0\}^h\oplus\mathcal{L}_\mathrm{scale}\subseteq\mathcal{H}$ contains all differences $(0;e_p-e_{p^\prime})$ with $p,p^\prime\in U_\alpha$, and for every $p\in W_\alpha$ there is an $i\in\mathcal{J}\cap C_\alpha$ such that $p=\iota(i)$, and $(0;e_p-e_{\sigma(i)})\in\mathcal{L}_\mathrm{IC}$ with $\sigma(i)\in U_\alpha$. These differences connect all vertices of $V_\alpha$, so
\begin{equation}
\{0\}^h\oplus\mathcal{L}(H_\mathrm{aug})\subseteq\mathcal{H}+\mathcal{L}_\mathrm{IC},
\end{equation}
and consequently,
\begin{equation}
(\mathcal{H}+\mathcal{L}_\mathrm{IC})\cap\ker(\chi_\mathrm{aug})
=\{0\}^h\oplus\mathcal{L}(H_\mathrm{aug}).
\end{equation}
On $\{0\}^h\oplus\mathbb{Z}^U$, $\chi_\mathrm{aug}$ factors through the component map of $H_\mathrm{scale}$, since every component of $H_\mathrm{scale}$ lies in a component of $H_\mathrm{aug}$. The image of $\mathcal{H}$ is therefore obtained by projecting $\Lambda_\mathrm{rate}$ onto the components of $H_\mathrm{aug}$. Noting that the initial-condition differences are annihilated by $\chi_\mathrm{aug}$, the quotient-image identity in the proof of \Cref{maintheorem}, together with \Cref{rassignedprop}, gives
\begin{equation}
\chi_\mathrm{aug}(\mathcal{H}+\mathcal{L}_\mathrm{IC})
=\pi_\sigma^\mathrm{aug}(\bar{\Lambda}_\mathrm{rate}).
\end{equation}
Applying \Cref{rankadditivitylemma} first to $\chi_\mathrm{aug}$ and then to $\pi_\kappa$ yields
\begin{equation}
\operatorname{rank}_\mathbb{Z}(\Lambda_\mathrm{IVP})
=h+s_\mathrm{tot}-c_\mathrm{aug}+r_\mathrm{IVP}
=N_\mathrm{IVP}(\sigma)+1.
\end{equation}

The homomorphism $f(v)=v+\deg_\kappa(v)v_T$ fixes $\mathcal{L}_\mathrm{IC}$ and sends $\Gamma$ onto $\tilde{\Gamma}$, implying that
\begin{equation}
f(\Lambda_\mathrm{IVP})=\tilde{\Gamma}_\mathrm{IVP}.
\end{equation}
By \Cref{lem:timedeg}, $\deg_\kappa(v_T)=-1$, with the kernel argument of \Cref{Gammalatticelemma} showing that the restriction of $f$ to $\Lambda_\mathrm{IVP}$ has rank-one kernel precisely when $v_T\in\Lambda_\mathrm{IVP}\otimes_\mathbb{Z}\mathbb{Q}$, and has trivial kernel otherwise, proving \eqref{eq:IVPtimecount}. If $v_T\in\Lambda_\mathrm{IVP}$, which every anchor satisfies, then $f(\Lambda_\mathrm{IVP})\subseteq\Lambda_\mathrm{IVP}\cap\ker(\deg_\kappa)$ and $f$ fixes every element of that intersection. Since $\deg_\kappa$ vanishes on $\mathcal{L}_\mathrm{IC}$, the intersection equals $\Delta\Gamma+\mathcal{L}_\mathrm{IC}$, proving \eqref{eq:IVPanchorlattice}. Finally, if $W_\alpha\subseteq U_\alpha$ for all $\alpha$, then $V_\alpha=U_\alpha$, $H_\mathrm{aug}=H_\mathrm{scale}$, $s_\mathrm{tot}=s$ and $r_\mathrm{IVP}=r_\mathrm{rate}$, so the two counts coincide.
\end{proof}
A time multiplier formed from a basis of $\Delta\Gamma+\mathcal{L}_\mathrm{IC}$ keeps an anchor valuation in $\Lambda_\mathrm{IVP}$, although it may leave $\Gamma_\mathbb{Q}$. Thus the broader time-scale condition is required when polynomialising coefficients in a full-IVP basis, and the exact joint lattice remains unchanged.\\~\\
We now let $H_0$ be the graph on $\iota(\mathcal{J})$ formed by the cliques on $W_\alpha$, and let $H_0^+$ have vertex set $U_\mathrm{cand}$ and the same edges, setting $c(H_0)=0$ when $\mathcal{J}=\emptyset$.
\begin{theorem}[A lower bound and its attainment for scale assignments]\label{thm:ivpmin}\label{proof:ivpmin}
Suppose that $r_\mathrm{IVP}(\sigma)$ is constant on a nonempty admissible family $\Sigma$, with common value $r_\mathrm{IVP}$, with this holding if $\bar{\Lambda}_\mathrm{rate}=\{0\}$, in particular if all rate symbols are distinct. Then for $\sigma\in\Sigma$,
\begin{equation}
N_\mathrm{IVP}(\sigma)\geq N_\mathrm{IVP}^\ast:=
h-1+|\iota(\mathcal{J})|-c(H_0)+r_\mathrm{IVP}, \label{eq:NIVPmin}
\end{equation}
with equality holding when every $V_\alpha(\sigma)$ lies in one component of $H_0^+$. If $\Sigma$ contains such an assignment, these are precisely its parameter-minimal assignments, and the minimum is $N_{\mathrm{IVP}}^\ast$; otherwise, the constrained minimum is strictly larger.
\end{theorem}
\begin{proof}
We treat the two graph lattices as sublattices of $\mathbb{Z}^{U_\mathrm{cand}}$, adjoining unused vertices as isolated vertices where necessary. Since $W_\alpha\subseteq V_\alpha$ for every $\alpha$,
\begin{equation}
\mathcal{L}(H_0^+)\subseteq\mathcal{L}(H_\mathrm{aug}(\sigma)).
\end{equation}
Noting that adjoining isolated vertices does not change a graph lattice or its rank, \Cref{graphranklemma} gives
\begin{equation}
|\iota(\mathcal{J})|-c(H_0)
\leq s_\mathrm{tot}(\sigma)-c_\mathrm{aug}(\sigma),
\end{equation}
and \eqref{eq:NIVPmin} follows from \eqref{eq:NIVP}. The lattice $\mathcal{L}(H_0^+)$ is the kernel of the component-sum homomorphism into a free abelian group, and is therefore saturated. Thus, under the displayed inclusion, equality of ranks consequently holds if and only if the two graph lattices coincide. By \Cref{graphranklemma}, this occurs if and only if $e_p-e_{p^\prime}\in\mathcal{L}(H_0^+)$ for every $p,p^\prime\in V_\alpha$ and every $\alpha$, which is equivalent to the stated component condition. Additionally, the final assertions follow because $\Sigma$ is a nonempty subset of a finite set of assignments.
\end{proof}
For a component with $W_\alpha\neq\emptyset$, the equality condition requires its species scales to lie in the $H_0$-component containing $W_\alpha$. If $W_\alpha=\emptyset$, they must all lie in one component of $H_0^+$, possibly an isolated candidate symbol.
\begin{corollary}\label{cor:default}
Under the hypotheses of \Cref{thm:ivpmin}, suppose that $\sigma\in\Sigma$ satisfies $\sigma(i)=\iota(i)$ for $i\in\mathcal{J}$, uses symbols from $W_\alpha$ for the remaining species of each component with $W_\alpha\neq\emptyset$, and is constant on components with $W_\alpha=\emptyset$. Then all positive initial values are normalised to one and $N_\mathrm{IVP}(\sigma)=N_\mathrm{coef}(\sigma)=N_\mathrm{IVP}^\ast$.
\end{corollary}
\begin{proof}
Each $V_\alpha$ is $W_\alpha$ or a singleton, and $W_\alpha\subseteq U_\alpha$, so \Cref{thm:ivpmin,thm:ivpcount} apply.
\end{proof}
\begin{numberedremark}[Coefficient-only minimisation]\label{rem:coefonly}
It is important to note that the smaller coefficient count in \Cref{mingentheorem} need not reduce the full-IVP count. Furthermore, the shared-rate example in \Cref{ex:rrate-ivp} shows how changes in the augmented graph and rate-rank term can compensate when $r_\mathrm{IVP}$ is not constant.
\end{numberedremark}

\begin{example}[\Cref{ex:rrate} continued]\label{ex:rrate-ivp}
Letting $X_1(0)=a$ and $X_2(0)=b$, we see that the graph $H_0$ consists of two isolated vertices. For $\sigma=(a,b)$, the augmented scale graph has two components and $r_\mathrm{IVP}=1$, so $N_\mathrm{IVP}=3+2-2+1-1=3$. For $\sigma=(a,a)$, the initial-condition ratio $b/a$ adds the edge joining $a$ and $b$, so $c_\mathrm{aug}=1$, $r_\mathrm{IVP}=0$, and again $N_\mathrm{IVP}=3$. With $K_T=\kappa_1^{-1}$, a Laurent basis in the latter case is $\{\kappa/\kappa_1,\kappa_2/\kappa_1,b/a\}$, with $r_\mathrm{IVP}$ not being constant across the two assignments, so \Cref{thm:ivpmin} does not apply. We see that the increase in the graph contribution to \eqref{eq:NIVP} is compensated by the decrease in the rate-lattice contribution.
\end{example}
\subsection{Comparison with the Buckingham Count}\label{subsec:buckingham}
The Buckingham count for these symbols, equal to the rank of their dimensionless monomial lattice, is $N_\Pi(\sigma):=h+s_\mathrm{tot}-\operatorname{rank}(\mathcal{D}_\sigma)$, where $\mathcal{D}_\sigma$ is the restriction of $\mathcal{D}$ to $\mathbb{Z}^h\oplus\mathbb{Z}^{U_\mathrm{tot}}$ \cite{Buckingham1914}.
\begin{proposition}\label{prop:buckingham}
For every admissible $\sigma$, $N_\mathrm{IVP}(\sigma)\leq N_\Pi(\sigma)$, with equality holding if $G$ is connected and all species and symbols in $U_\mathrm{tot}$ have a common nonzero dimension.
\end{proposition}
\begin{proof}
Under anchor normalisation, all coefficient and initial-value valuations belong to $\ker(\mathcal{D}_\sigma)$, proving the inequality by rank comparison. With a common nonzero dimension $\delta$, the rate dimensions lie in the span of $e_T$ and $\delta$ and have time component $-1$, while scale symbols have dimension $\delta$ and zero time component. Hence, $\operatorname{rank}(\mathcal{D}_\sigma)=2$, and if $G$ is also connected, then $c_\mathrm{aug}=1$. For reactions $j,k$ sharing a rate symbol, subtracting \eqref{eq:dimconsistency} gives
\begin{equation*}
\left(\sum_{\ell=1}^n\nu_{\ell j}-\sum_{\ell=1}^n\nu_{\ell k}\right)\delta=0.
\end{equation*}
Since $\delta\neq 0$, their total reactant orders coincide, so their reaction-scale vectors agree after projection to the single augmented component, implying that $r_\mathrm{IVP}=0$, and $N_\mathrm{IVP}=N_\Pi=h+s_\mathrm{tot}-2$.
\end{proof}
We note that for a common nonzero concentration dimension, the same calculation gives
\begin{equation}
N_\Pi-N_\mathrm{IVP}=c_\mathrm{aug}-r_\mathrm{IVP}-1,
\end{equation}
so that disconnectedness alone does not imply that $N_\Pi>N_\mathrm{IVP}$.
\begin{numberedremark}[Scope of minimality]\label{rem:scope}
These counts concern Laurent-monomial representations of fixed reaction-channel contributions and initial values, with species scales drawn from the prescribed symbols. Combining channels, imposing parameter relations, eliminating conserved quantities or permitting more general changes of variables may alter the count, but these possibilities are beyond the scope of this work.
\end{numberedremark}

\section{Algorithmic Implementation with the CRNNondim Python Package\label{algosection}}
\subsection{CRN Structural Preprocessing}
Given a CRN in the form of \eqref{maincrnform}, it is prudent to preprocess the CRN to extract the participation matrix, $P$, the species co-participation graph $G$, its components, the lattice $\bar{\Lambda}_\mathrm{rate}$, the initial-value graph $H_0$, and, for a fixed scale assignment, the exponent matrix $A$. From this, we can evaluate the admissible scale assignments and, on any family on which $r_\mathrm{IVP}$ is constant and the lower bound in \Cref{thm:ivpmin} is attainable, minimise the number of parameters of the full IVP.  In the general case, we aim to minimise \eqref{eq:NIVP} over that finite family, with the coefficient-only objective remaining $s-c_\mathrm{scale}+r_\mathrm{rate}$, as in \Cref{maintheorem,mingentheorem}. Additionally, we can construct the exponent matrix, $A$, and exponent lattice, $\Gamma$, corresponding to the monomial coefficients, $c_{ij}$, in terms of the distinct rate constants and used scale symbols. An algorithm for this structural preprocessing in the coefficient-only case is provided in Algorithm C.1.
\subsection{Scale Selection}
An important aspect of non-dimensionalisation is to keep all state variables as close to unity as possible over a relevant time window, which can improve conditioning in numerical ODE integration and related computations \cite{Hairer1996}. Intuitively, a dimensionless variable with magnitude $10$ should be penalised just as much as one with magnitude $0.1$, with this reciprocity naturally motivating objectives of the type $\max\limits_{i\in \mathcal{I}}\left|\ln(x_i(T/K_T))\right|$, where the base of the logarithm can be chosen to be any number greater than 1. We assume numerical values for the rate constants and initial conditions, and a solution with $X_i(T)>0$ for every species throughout $[T_\mathrm{min},T_\mathrm{max}]$, where $0\leq T_\mathrm{min}\leq T_\mathrm{max}$. All extrema below are taken over this interval, and throughout this section, logarithms of dimensional quantities denote logarithms of their positive numerical values in a fixed, consistent system of units. Writing $X_i^{\min}=\min_T X_i(T)$ and $X_i^{\max}=\max_T X_i(T)$, we set
\begin{equation}
\alpha_i:=\ln X_i^{\min},\qquad \beta_i:=\ln X_i^{\max},\qquad
m_i:=\tfrac12(\alpha_i+\beta_i),\qquad h_i:=\tfrac12(\beta_i-\alpha_i).
\label{definetempvars}
\end{equation}
For $L_\ell:=\ln a_\ell$ and $K_i=a_{\sigma(i)}$, we define
\begin{align}
\Phi_i(T;\sigma)&:=\ln\!\left(\frac{X_i(T)}{K_i}\right)=\ln X_i(T)-L_{\sigma(i)},\label{Phidef}\\
r_i(\sigma)&:=\max_T|\Phi_i(T;\sigma)|=h_i+|L_{\sigma(i)}-m_i|,\label{maxdef}
\end{align}
where the last identity follows since the maximum absolute deviation over a real interval is attained at an endpoint, which we prove in Lemma B.1 and Corollary B.1. Thus, for fixed trajectories $X_i(T)$ and candidate scales $(a_\ell)$, minimising the objective $R(\sigma):=\max_{i\in\mathcal{I}}r_i(\sigma)$ is exactly the same as minimising the worst log-deviation over all species, leading to \Cref{prop:supnorm-exact}, which we prove in Section B.2.
\begin{proposition}\label{prop:supnorm-exact}
Assume that the allowable set of species scales per species satisfies $\emptyset\neq S_i\subseteq U_\mathrm{cand}$ for all $i\in\mathcal{I}$. For each species $i$, define the best achievable per-species error
  \begin{equation}
  r_i^\ast := \min_{\ell\in S_i} \left(h_i + |L_\ell - m_i|\right),
\end{equation}
and set 
\begin{equation}
  R^\ast := \max_{i\in\mathcal{I}} r_i^\ast.
\end{equation}
Then,
\begin{enumerate}
  \item[(a)] $R^\ast$ is the global minimum of $R(\sigma)$ over all admissible scale assignments:
  \begin{equation}
    R^\ast = \min_{\sigma} R(\sigma),
  \end{equation}
  where $\sigma$ ranges over all maps $\sigma:\mathcal{I}\to U_\mathrm{cand}$ with $\sigma(i)\in S_i$ for every $i$.
  \item[(b)] For every assignment $\sigma$ with $R(\sigma)=R^\ast$, $r_i(\sigma)\geq r_i^\ast$ for all $i$, and there exists at least one index $i$ with $r_i(\sigma)=R^\ast$.
  \item[(c)] For each $i\in\mathcal{I}$, the set
  \begin{equation}
  \mathcal{S}_i^\ast := \set{\ell\in S_i \mid h_i + |L_\ell - m_i| \leq R^\ast}
  \end{equation}
  is nonempty, and every assignment $\sigma$ with $\sigma(i)\in\mathcal{S}_i^\ast$ for all $i\in\mathcal{I}$ satisfies $R(\sigma)\leq R^\ast$.
\end{enumerate}
In particular, computing $r_i^\ast$ and $R^\ast$ first, and then restricting to assignments with $\sigma(i)\in\mathcal{S}_i^\ast$, attains the global minimum of the sup-norm objective $R(\sigma)$.
\end{proposition}
Given the admissible species scale sets $S_i$ and numerical ranges of $X_i(T)$ in the time window $T\in[T_\mathrm{min},T_\mathrm{max}]$, Algorithm C.2 implements the exact discrete optimisation algorithm to characterise all assignments that minimise $R(\sigma)$ over all admissible scale assignments.\\~\\
In general, there are many assignments with $\sigma(i)\in\mathcal{S}_i^\ast$ for all $i\in\mathcal{I}$. While each such $\sigma$ attains the same maximum logarithmic deviation $R^\ast$, they can significantly differ in how the per-species deviations $r_i(\sigma)$ are distributed. As such, it is natural to find the optimal assignment via lexicographic refinement. Specifically, for every admissible $\sigma$, we sort $(r_i(\sigma))_{i\in\mathcal{I}}$ in non-increasing order and select those assignments that minimise this vector in lexicographic order, thereby minimising, successively, the second-largest deviation, the third-largest, and so on.\\~\\
We define $\Sigma^\ast:= \set{\sigma \mid R(\sigma) = R^\ast}$. It is important to note that by \Cref{prop:supnorm-exact}, $\Sigma^\ast=\mathcal{S}_1^\ast\times\cdots\times\mathcal{S}_n^\ast$, so it is represented by the nonempty label sets returned by Algorithm C.2. For each $\sigma\in\Sigma^\ast$, we define $\bm{r}(\sigma) := (r_i(\sigma))_{i\in\mathcal{I}}\in\mathbb{R}_{\geq 0}^n$, and let $\bm{r}^\downarrow(\sigma)\in\mathbb{R}_{\geq 0}^n$ be the vector obtained by sorting the entries of $\bm{r}(\sigma)$ in non-increasing order. We equip $\mathbb{R}_{\geq 0}^n$ with the usual lexicographic order, so that for $\bm{u},\bm{v}\in\mathbb{R}_{\geq 0}^n$, we write $\bm{u}\prec_{\mathrm{lex}}\bm{v}$ if at the first index $k$ where $u_k\neq v_k$, we have $u_k < v_k$.
\begin{proposition}\label{prop:lexi-exact}
Let $\Sigma^\ast$ and $\bm{r}^\downarrow(\sigma)$ be defined as above. Then, 
\begin{enumerate}
  \item[(a)] There exists at least one assignment $\sigma^{\mathrm{lex}}\in\Sigma^\ast$ such that
  \begin{equation}
  \bm{r}^\downarrow(\sigma^{\mathrm{lex}}) \preceq_{\mathrm{lex}} \bm{r}^\downarrow(\sigma) \quad\text{for all }\sigma\in\Sigma^\ast.
  \end{equation}
  \item[(b)] Any such $\sigma^{\mathrm{lex}}$ is a global minimiser of the lexicographic objective $\bm{r}^\downarrow(\sigma)$ over all admissible assignments.
\end{enumerate}
\end{proposition}
Given the admissible sets $\{\mathcal{S}_i^\ast\}$ produced by Algorithm C.2, Algorithm C.3 implements a lexicographic refinement that determines $\sigma^{\mathrm{lex}}\in\Sigma^\ast$, as defined in \Cref{prop:lexi-exact}, which we prove in Section B.2. We note that for a Cartesian family, all per-species minima can be attained simultaneously. Thus, writing $\displaystyle A_i^\ast:=\operatorname*{argmin}_{\ell\in\mathcal{S}_i^\ast}(h_i+|L_\ell-m_i|)$, the full lexicographically minimal family is $\Sigma^{\mathrm{lex}}=\prod\limits_{i\in\mathcal{I}}A_i^\ast$, represented without enumerating this product. Distinct scale labels remain distinct even when their numerical values agree, so every such tied label is retained for the subsequent symbolic search.\\~\\
In particular, \Cref{prop:supnorm-exact} and \Cref{prop:lexi-exact} together show that Algorithm C.2 and Algorithm C.3 exactly minimise $\max\limits_{i\in\mathcal{I}}\max\limits_{T\in[T_{\mathrm{min}},T_{\mathrm{max}}]} \left|\ln(x_i(T/K_T))\right|$, and hence $\max\limits_{i\in\mathcal{I}}\max\limits_{T\in[T_{\mathrm{min}},T_{\mathrm{max}}]} \left|\log_{10}(x_i(T/K_T))\right|$ over the discrete set of admissible assignments $\sigma$. We note that the extrema $X_i^{\min}$ and $X_i^{\mathrm{max}}$ are numerically obtained from ODE integrators, so that optimality is exact relative to these approximations and the comparison tolerances.\\~\\
Alternatively, when positive initial concentrations fix their species scales, one may wish to select the unfixed species scales such that species in $\mathcal{I}^\ast:=\{i:X_i(0)=0\}$ may be scaled so that $\max_Tx_i=X_i^{\max}/K_i$ does not exceed unity while remaining as close to it as possible. Thus, for each $i\in\mathcal{I}^\ast$, we define
\begin{equation}
\mathcal{S}_i^{(1)} := \begin{cases}
\operatorname*{argmin}\limits_{\ell\in S_i \mid a_\ell \geq X_i^{\mathrm{max}}} a_\ell,
& \text{if } \exists \ell\in S_i \text{ such that } a_\ell\geq X_i^{\mathrm{max}},\\[2mm]
\operatorname*{argmax}\limits_{\ell\in S_i} a_\ell,
& \text{otherwise}.
\end{cases}
\end{equation}
The second case arises when none of the allowable scales is sufficiently large to ensure that the corresponding non-dimensionalised trajectory does not exceed unity, in which case choosing the largest available scale instead makes the value of $X_i^{\mathrm{max}}/K_i>1$ as small, and hence as close to unity, as possible. As with the deviation-based criterion above, distinct scale labels whose numerical values coincide remain tied and are retained for the subsequent symbolic optimisation, since, although they are indistinguishable under this numerical scale-selection criterion, they need not be equivalent in the resulting symbolic non-dimensionalisation.
\subsection{Basis simplification}\label{subsec:unimod-simplification}
From \Cref{unimodularmain}, the choice of $\mathbb{Z}$-basis $\{g_1,\dots,g_{d-1}\}$ for the lattice $\tilde{\Gamma}\cong\mathbb{Z}^{d-1}$ is not unique. Given an anchor-normalised time scaling $K_T\in \mathcal{M}$, any unimodular transformation $U \in \operatorname{GL}_{d-1}(\mathbb{Z})$ yields a new basis. Let $E:=\mathcal{P}$ and let $g_1,\ldots,g_N$ be a basis of the relevant coefficient or IVP lattice, where $N=d-1$ in the coefficient-only case, and let $B$ contain its valuation columns and let $V$ contain the coefficient-coordinate rows. By \Cref{unimodularmain,laurentbasisexistence}, a unimodular transformation gives
\begin{equation}
g_k^\prime=\prod_{i=1}^N g_i^{(U^{-1})_{ki}},\qquad
V^{(U)}=VU,\qquad B(U)=BU^{-\intercal},\qquad U\in\mathrm{GL}_N(\mathbb{Z}).
\end{equation}
It is thus natural to choose the unimodular matrix $U$ so that the exponent vectors $V^{(U)}_e$ are as sparse and small as possible, since this both simplifies subsequent analysis and allows us to automatically generate more canonical choices of generators. This motivates an $\ell^1$-type penalty
\begin{equation}
F(U):=\sum_{e\in E}\|V_e^{(U)}\|_1 =\sum_{e\in E}\sum_{k=1}^N|(VU)_{e,k}|.\label{eq:L1-cost}
\end{equation}
Clearly, a small value for $F(U)$ implies that for most of the coefficients $\tilde{c}_e$, only a few generators $g^\prime_k$ appear, with small integer exponents. Although $F$ is convex on the ambient real matrix space, $\mathrm{GL}_N(\mathbb{Z})$ is nonconvex and infinite for $N \geq 2$, and decision versions of related sparsification problems are NP-complete \cite{osti_5633555}. We therefore use a greedy search over unit shears $I_N+\xi E_{pq}$, where $E_{pq}$ is the matrix unit, $p\neq q$, and $\xi=\pm1$. These correspond to adding $\xi\operatorname{col}_p(V)$ to $\operatorname{col}_q(V)$ and replacing the generator $g_p$ by $g^\prime_p = g_p g_q^{-\xi}$ while leaving the other generators unchanged.\\~\\
We present an algorithmic description of a greedy local search to minimise \eqref{eq:L1-cost} in Algorithm C.4. We note that unlike Algorithms C.2 and B.3, Algorithm C.4 does not guarantee a globally minimal value of $F(U)$ over all $U\in\mathrm{GL}_{N}(\mathbb{Z})$, but instead finds a basis that is locally optimal under the restricted set of elementary unimodular moves considered. In practice, this is sufficient to find bases that generate exponent vectors for the $\tilde{c}_{ij}$ that are significantly sparser and smaller than those produced by just considering the SNF, resulting in much simpler monomial coefficients. Henceforth, statements that a basis is as sparse and small as possible are understood relative to this restricted set of moves.
\subsection{Converting between Laurent and Polynomial Bases}\label{subsec:polynomial-simplification}
The Laurent-to-polynomial conversion in \Cref{thm:poly-shift} serves as the foundation for the conversion used in practice to generate polynomial bases of $\{\tilde{c}_{ij}\}$, with an algorithmic description for the global coordinate-wise shift being described in Algorithm C.5. It is important to note that minimising the Laurent cost before this shift need not minimise the final polynomial cost, therefore motivating evaluating every candidate after its required time-scale transformation. With $a_T=\operatorname{val}(K_T)$ for an initial anchor, define
\begin{align}
R(U)&:=VU,& B(U)&:=BU^{-\intercal},& M_k(U)&:=\min_{e\in E}R(U)_{e,k},\label{eq:poly-candidate}\\
Q(U)&:=R(U)-\mathbf{1}M(U),& a_T(U)&:=a_T-B(U)M(U)^\intercal,\label{eq:poly-time-valuation}
\end{align}
where $M(U) = (M_1(U), \dots, M_N(U))$ is a row vector and $\mathbf{1}$ is the all-ones vector in $\mathbb{R}^{|E|}$. In particular, $Q(U)$ is non-negative, and $a_T(U)$ is the valuation of $K_T\prod_{k=1}^{N}(g^\prime_k)^{-M_k(U)}$. Then, the corresponding polynomial penalty is
   \begin{equation}
\mathcal{F}_\mathrm{poly}(U):=\left(\sum_{e,k}Q(U)_{e,k},\ \sum_{p,k}|B(U)_{p,k}|,\ \|a_T(U)\|_1\right),\label{eq:poly-score}
\end{equation}
ordered lexicographically as $(F_\mathrm{mono},F_\mathrm{gen},F_\mathrm{time})$. For a full-IVP basis, the initial-value rows $H$ become $HU$ and are not shifted, with the zero anchor supplying the shifted row $-M(U)$, whose differences with the other shifted coefficient rows recover $VU$. Together with $HU$, these rows therefore generate the same IVP lattice, and, if $N_\mathrm{IVP}>N_\mathrm{coef}$, the coefficient-only rank can increase, while the time valuation remains in $(\Gamma+\mathcal{L}_\mathrm{IC})\otimes\mathbb{Q}$. \\~\\
For each anchor, the best initial transformation under \eqref{eq:poly-score} is selected from the sparsity-first, coefficient-weighted $\ell^1$, and identity seeds, retaining this order for ties, and a single greedy descent then tests sign flips followed by signed unit shears. A move is accepted only if the complete post-shift $(F_{\mathrm{mono}},F_{\mathrm{gen}},F_{\mathrm{time}})$ strictly lexicographically decreases, with this procedure being described in Algorithm C.6. Since $\mathbb{N}^3$ admits no infinite strictly descending lexicographic chain, the descent terminates, so despite exact selection among evaluated anchors and tied assignments, the local basis search does not certify a global minimum over all unimodular bases.

\subsection{Processing a CRN under a Fixed Species Non-dimensionalisation}\label{subsec:fixednondim}
Given a specific choice of non-dimensionalisation, we aim to symbolically build the coefficient structure in \eqref{setcijeqn} for that assignment $\sigma$. We first note that it is more convenient to represent the monomials $c_{ij}$ from \eqref{setcijeqn} in a log-linear form. For a fixed $\sigma$, we define
\begin{equation}
\begin{aligned}
\ln c_{ij}&=\ln\kappa_{\rho(j)}+C_{(i,j)}x_\mathrm{sol},\\
C_{(i,j),\ell}&:=\nu_{\ell j}-\delta_{i\ell},\qquad
x_\mathrm{sol}:=(\ln K_i)_i,
\end{aligned}
\end{equation}
where $\delta_{i\ell}$ is the Kronecker delta, so that the logarithm of each coefficient is affine in the log-scales. We encode the corresponding valuation rows $v^{(e)}$ in $V_\mathrm{abs}$, using the common ambient lattice $\mathbb{Z}^h\oplus\mathbb{Z}^{U_\mathrm{tot}}$, with zero entries for symbols absent from the coefficients. For an anchor $e_\ast$, form
\begin{equation}
D_e=v^{(e)}-v^{(e_\ast)},\qquad
(A_\mathrm{IC})_i=(0;e_{\iota(i)}-e_{\sigma(i)})\quad(i\in\mathcal{J}),\qquad
D_\mathrm{IVP}=\begin{bmatrix}D\\A_\mathrm{IC}\end{bmatrix},
\end{equation}
so that its rows generate $\Delta\Gamma+\mathcal{L}_\mathrm{IC}$. Computing the SNF
\begin{equation}
TD_\mathrm{IVP}Q=\operatorname{diag}(\alpha_1,\ldots,\alpha_N,0,\ldots,0),\qquad
N=\operatorname{rank}_{\mathbb{Z}}(D_\mathrm{IVP}),
\end{equation}
with unimodular $T,Q$ gives a row basis $B_\mathrm{IVP}$ from the first $N$ rows of $TD_\mathrm{IVP}$ and integer coordinates $W=(T^{-1})_{[:,1:N]}$ satisfying $D_\mathrm{IVP}=WB_\mathrm{IVP}$. Setting $B=B_\mathrm{IVP}^\intercal$ and splitting $W$ into coefficient and initial rows gives $D=VB^\intercal$ and $A_\mathrm{IC}=HB^\intercal$. Omitting $A_\mathrm{IC}$ recovers the coefficient-only construction with rank $d-1$, detailed in Section B.3.\\~\\
We can then simplify this exponent matrix via Algorithm C.4 or Algorithm C.6, depending on whether Laurent or polynomial representations of $\{\tilde{c}_{ij}\}$ are required, obtaining a unimodular matrix $U\in \mathrm{GL}_{N}(\mathbb{Z})$ and a simplified exponent matrix $V^{(U)}$. Changing the anchor to $e_0$ replaces $V$ by $V-\mathbf{1}V_{e_0}$ and sets $a_T=-(v^{(e_0)})^\intercal$, corresponding to $K_T=c_{e_0}^{-1}$. After the chosen basis transformation and any coefficient shift, the returned matrices satisfy
\begin{equation}
V^\prime(B^\prime)^\intercal=V_\mathrm{abs}+\mathbf{1}(a_T^\prime)^\intercal,\qquad
H^\prime(B^\prime)^\intercal=A_\mathrm{IC}.
\end{equation}
Each choice of anchor produces one family, with its own $K_T$, generators, and coefficient representation. After unimodular simplification and, if desired, polynomialisation, we compute the total monomial $\ell^1$ score, $F_\mathrm{mono}$, and consider only those families that attain the minimum values. We then apply a similar filtering based on the total generator $\ell^1$ score, $F_\mathrm{gen}$, and then the $K_T$ expression $\ell^1$ score, $F_\mathrm{time}$, and retain only those families, including their induced $K_T$, that minimise these scores. An algorithmic description of the full procedure in the coefficient-only case is provided in Algorithm C.7.
\subsection{Determining an Optimal Non-dimensionalisation for a CRN}
We combine Algorithms C.1 to B.7 into a single algorithm, shown in Algorithm C.8, which selects a lexicographically minimal family among those constructed for a given CRN with numerical rate constants and initial conditions. This choice of non-dimensionalisation keeps the non-dimensionalised variables as close to unity as possible among the structurally minimal assignments, and subsequently filters the generated families by the monomial, generator, and $K_T$ expression $\ell^1$ scores. Additionally, to break ties between different families, we choose the family for which $t_\mathrm{max} := T_\mathrm{max}/K_T$ is as close to $1$ as possible, formally minimising $T_\mathrm{cost}:=\left|\log_{10}(T_\mathrm{max}/K_T)\right|$. 

\subsection{The CRNNondim Python Package}
We have implemented the coefficient-based versions of the algorithms in Algorithms C.1 to B.8 in the Python package CRNNondim, available on PyPI and at \url{https://github.com/georgio-hawi/crnnondim}. The implementation uses SymPy for symbolic computer algebra \cite{Meurer2017}, NumPy for numerical linear algebra \cite{Harris2020}, and SciPy's Radau integrator for numerical integration of CRNs \cite{Virtanen2020}, an implicit Runge--Kutta method well suited to stiff systems such as CRNs \cite{Hairer1996,Hairer1999}. The IVP lattice augmentation represents a mathematical extension whose implementation falls outside the scope of this work.
\section{Examples\label{examplessection}}
In this section, we provide worked examples that illustrate how CRNNondim applies the theory and algorithms developed in \Cref{resultssection} and \Cref{algosection} to non-dimensionalise several well-known biological CRNs, with detailed exponent matrices and numerical details provided in Section D. In each example, direct evaluation gives $\bar{\Lambda}_\mathrm{rate}=\{0\}$, so that $r_\mathrm{rate}=r_\mathrm{IVP}=0$.
\numberwithin{equation}{subsection}
\subsection{Reversible Michaelis-Menten Scheme}\label{subsec:mm}
The reversible Michaelis-Menten scheme is an archetypal CRN that underpins modern enzyme kinetics and metabolic modelling, and is ubiquitous as a minimal motif for enzyme-catalysed reactions throughout biochemistry \cite{michaelis1913kinetik, Johnson2011, CornishBowden2015}. It is described by the reaction scheme
\begin{equation}
\mathrm{E} + \mathrm{S} \underset{k_{-1}}{\stackrel{k_{1}}{\rightleftharpoons}} \mathrm{C} \underset{k_{-2}}{\stackrel{k_{2}}{\rightleftharpoons}} \mathrm{E} + \mathrm{P}, \label{michaelis-menteneqns}
\end{equation}
where $\mathrm{S}$, $\mathrm{E}$, $\mathrm{C}$, and $\mathrm{P}$ correspond to a substrate, enzyme, enzyme-substrate complex, and product, respectively. Denoting $\mathbf{X}:= [[\mathrm{S}]$, $[\mathrm{E}]$, $[\mathrm{C}]$, $[\mathrm{P}]]^\intercal$ $\in \mathbb{R}^{4}_{ \geq 0}$ as the vector of species concentrations, we can express \eqref{michaelis-menteneqns} as a system of ODEs.
We have that by the law of mass action, \eqref{michaelis-menteneqns} is expressible as 
\begin{equation}
\begin{aligned}
\frac{dX_{1}}{dT} &= \fitRHS{-k_{1}X_{1}X_{2} + k_{-1}X_{3}}, \\
\frac{dX_{2}}{dT} &= \fitRHS{-k_{1}X_{1}X_{2} + k_{-1}X_{3} + k_{2}X_{3} -k_{-2}X_{2}X_{4}}, \\
\frac{dX_{3}}{dT} &= \fitRHS{k_{1}X_{1}X_{2} -k_{-1}X_{3} -k_{2}X_{3} + k_{-2}X_{2}X_{4}}, \\
\frac{dX_{4}}{dT} &= \fitRHS{k_{2}X_{3} -k_{-2}X_{2}X_{4}},
\end{aligned}
\quad\implies\quad
\begin{aligned}
\frac{d\mathbf{X}}{dT}
&= \underbrace{\begin{bmatrix}
-1 & 1 & 0 & 0 \\
-1 & 1 & 1 & -1 \\
1 & -1 & -1 & 1 \\
0 & 0 & 1 & -1
\end{bmatrix}}_{S}
\underbrace{\begin{bmatrix}
k_{1}X_{1}X_{2} \\
k_{-1}X_{3} \\
k_{2}X_{3} \\
k_{-2}X_{2}X_{4}
\end{bmatrix}}_{V(\mathbf{X})}.
\end{aligned}
\label{eq:michaelismentendim}
\end{equation}

We take $\mathbf{X}(0)=(s_0,e_0,0,0)^\intercal$, with $U_\mathrm{cand}=\{s_0,e_0\}$. Since the four reaction-rate symbols are distinct and $G$ is connected, so that $H_0$ consists of the edge joining $s_0$ and $e_0$, \Cref{thm:ivpmin} gives $N^\ast_\mathrm{IVP}=4-1+2-1=4$, equal to the Buckingham count $6-2$. With $K_1=s_0$ and $K_2=e_0$, every assignment of the remaining scales from $U_\mathrm{cand}$ has four coefficient generators, while choosing $K_i=s_0$ for all $i$ instead reduces the coefficient count to three but introduces $x_2(0)=e_0/s_0$, leaving the full IVP count equal to four.\\~\\
We first consider $\mathcal{K}=\operatorname{diag}(s_0,e_0,e_0,e_0)$, with anchor $(i_0,j_0)=(2,2)$ and $K_T=k_{-1}^{-1}$, the reference timescale for $\mathrm{C}\xrightarrow{k_{-1}}\mathrm{E}+\mathrm{S}$. Using \Cref{unimodularmain} and Algorithm C.4, we select a $\mathbb{Z}$-basis of $\Delta\Gamma$ that simplifies the coefficient exponents, giving the Laurent basis $\mathcal{G}=\{g_1,\dots,g_4\}$, where
\begin{equation}
g_1 = \frac{k_1 s_0}{k_{-1}}, \quad  g_2 = \frac{k_2}{k_{-1}}, \quad  g_3 = \frac{k_{-2}e_0}{k_{-1}}, \quad g_4 = \frac{e_0}{s_0},
\end{equation}
so that the non-dimensionalised reversible Michaelis-Menten scheme is expressible as
\begin{equation}
    \frac{d\mathbf{x}}{dt} =
    \underbrace{\begin{bmatrix}
-g_4 &  g_4 &  0 &  0 \\
-1   &  1   &  1 & -1 \\
 1   & -1   & -1 &  1 \\
 0   &  0   &  1 & -1
\end{bmatrix}}_{\tilde{S}}
\underbrace{\begin{bmatrix}
g_1 x_1 x_2 \\
x_3 \\
g_2 x_3 \\
g_3 x_2 x_4
\end{bmatrix}}_{\tilde{V}(\mathbf{x})}, \label{eq:mmnondim1}
\end{equation}
with $\mathbf{x}(0)=(1,1,0,0)^\intercal$. It can be calculated that for this configuration, $F_\mathrm{mono} = 11$, $F_\mathrm{gen}=10$, and $F_\mathrm{time}=1$, and all monomial coefficients are polynomial in the generators so that $\mathcal{G}$ is also a polynomial basis. This choice of generators has a clear physical interpretation, with the first three generators comparing the forward binding, catalytic and reverse binding fluxes, respectively, with the dissociation flux $k_{-1}X_3$, evaluated at $X_i=K_i$, while $g_4$ measures the ratio of initial enzyme to substrate concentrations. For the corresponding irreversible scheme, with $k_{-2}=0$, the classical quasi-steady-state parameter, with $K_M=(k_{-1}+k_2)/k_1$, is
\begin{equation}
    \frac{e_0}{s_0+K_M}=\frac{g_1g_4}{g_1+1+g_2},
\end{equation}
which is a rational function of the monomial generators.\\~\\
We now consider the assignment $\mathcal{K}=\operatorname{diag}(s_0,e_0,e_0,s_0)$, retaining the same anchor. Using \Cref{unimodularmain} and Algorithm C.4, we select a $\mathbb{Z}$-basis of $\Delta\Gamma$ that simplifies the coefficient exponents, giving the Laurent basis $\mathcal{G}=\{g_1,\dots,g_4\}$, with $g_1$, $g_2$ and $g_4$ unchanged and $g_3=k_{-2}s_0/k_{-1}$. The last row of $\tilde{S}$ becomes $(0,0,g_4,-g_4)$, giving
\begin{equation}
    \frac{d\mathbf{x}}{dt} =
    \underbrace{\begin{bmatrix}
 -g_4 &  g_4 & 0 & 0 \\
 -1   &  1   & 1 & -1 \\
  1   & -1   & -1& 1 \\
  0   &  0   & g_4 & -g_4
\end{bmatrix}}_{\tilde{S}}
\underbrace{\begin{bmatrix}
g_1 x_1 x_2 \\
x_3 \\
g_2 x_3 \\
g_3 x_2 x_4
\end{bmatrix}}_{\tilde{V}(\mathbf{x})}, \label{mmnondim2}
\end{equation}
with the same initial conditions $\mathbf{x}(0)=(1,1,0,0)^\intercal$. It can be calculated that for this configuration, $F_\mathrm{mono} = 13$, $F_\mathrm{gen}=10$, and $F_\mathrm{time}=1$, and all monomial coefficients are polynomial in the generators so that $\mathcal{G}$ is also a polynomial basis.
\subsection{SEIRV Model}\label{subsec:seirv}
The 5-dimensional SEIRV (susceptible-exposed-infectious-recovered-vaccinated) model provides a framework for analysing how vaccination influences infectious disease dynamics, extending upon the very widely used SIR model in modern mathematical epidemiology \cite{kermack1927contribution, Anderson1991, Hethcote2000, Poonia2022}. We denote $\mathrm{S}$, $\mathrm{E}$, $\mathrm{I}$, $\mathrm{R}$ and $\mathrm{V}$ as the susceptible, exposed (infected, not yet infectious), infectious, recovered (immune) and vaccinated (protected) individuals, respectively. With $\mathbf{X}=([\mathrm{S}],[\mathrm{E}],[\mathrm{I}],[\mathrm{R}],[\mathrm{V}])^\intercal$, the SEIRV model is expressible as
\begin{align}
\frac{dX_{1}}{dT} &= \fitRHS{\Lambda - \beta X_{1}X_{3} - \nu_V X_{1} + \omega X_{4} + \omega_V X_{5} - \mu X_{1}}, \label{SEIRVeq1}\\
\frac{dX_{2}}{dT} &= \fitRHS{\beta X_{1}X_{3} - \sigma X_{2} - \mu X_{2}}, \\
\frac{dX_{3}}{dT} &= \fitRHS{\sigma X_{2} - \gamma X_{3} - \mu X_{3}}, \\
\frac{dX_{4}}{dT} &= \fitRHS{\gamma X_{3} - \omega X_{4} - \mu X_{4}}, \\
\frac{dX_{5}}{dT} &= \fitRHS{\nu_V X_{1} - \omega_V X_{5} - \mu X_{5}}. \label{SEIRVeq5}
\end{align}
Here $\Lambda$ is the recruitment rate, $\mu$ the natural per-capita death rate, $\beta$ the transmission rate and $\nu_V$ the vaccination rate, while the rates $\sigma$ and $\gamma$ govern progression to infectiousness and recovery, while $\omega$ and $\omega_V$ describe loss of immunity and vaccine protection. The model can be interpreted as a CRN with reactions
\begin{equation}
\begin{aligned}
\mathrm{I}+\mathrm{S} &\xrightarrow{\beta} \mathrm{E}+\mathrm{I}, &
\emptyset &\xrightarrow{\Lambda} \mathrm{S}, &
\mathrm{R} &\xrightarrow{\omega} \mathrm{S}, &
\mathrm{S} &\xrightarrow{\mu} \emptyset, &
\mathrm{S} &\xrightarrow{\nu_V} \mathrm{V}, &
\mathrm{V} &\xrightarrow{\omega_V} \mathrm{S}, \\
\mathrm{E} &\xrightarrow{\mu} \emptyset, &
\mathrm{E} &\xrightarrow{\sigma} \mathrm{I}, &
\mathrm{I} &\xrightarrow{\gamma} \mathrm{R}, &
\mathrm{I} &\xrightarrow{\mu} \emptyset, &
\mathrm{R} &\xrightarrow{\mu} \emptyset, &
\mathrm{V} &\xrightarrow{\mu} \emptyset,
\end{aligned}
\label{eq:seirv-reactions}
\end{equation}
so that \eqref{SEIRVeq1}--\eqref{SEIRVeq5} can equivalently be written as
\begin{equation}
\frac{d\mathbf{X}}{dT} =
    \underbrace{\begin{bmatrix}
-1 & 1 & 1 & -1 & -1 & 1 & 0 & 0 & 0 & 0 & 0 & 0 \\
1 & 0 & 0 & 0 & 0 & 0 & -1 & -1 & 0 & 0 & 0 & 0 \\
0 & 0 & 0 & 0 & 0 & 0 & 0 & 1 & -1 & -1 & 0 & 0 \\
0 & 0 & -1 & 0 & 0 & 0 & 0 & 0 & 1 & 0 & -1 & 0 \\
0 & 0 & 0 & 0 & 1 & -1 & 0 & 0 & 0 & 0 & 0 & -1
\end{bmatrix}}_{S}
\underbrace{\begin{bmatrix}
\beta X_{1}X_{3} \\
\Lambda \\
\omega X_{4} \\
\mu X_{1} \\
\nu_V X_{1} \\
\omega_V X_{5} \\
\mu X_{2} \\
\sigma X_{2} \\
\gamma X_{3} \\
\mu X_{3} \\
\mu X_{4} \\
\mu X_{5}
\end{bmatrix}}_{V(\mathbf{X})}.
\label{eq:seirvdim}
\end{equation}
We consider initial conditions of $\mathbf{X}(0)=(S_0,0,I_0,0,0)^\intercal$ and $U_\mathrm{cand}=\{S_0,I_0\}$, so that $h=8$. The graph $G$ is connected and $\mu$ is shared only by first-order death reactions whose reactant is their only active species, giving $\bar{\Lambda}_\mathrm{rate}=\{0\}$. The graph $H_0$ has one edge, so that $N^\ast_\mathrm{IVP}=8-1+2-1=8$, equal to the Buckingham count $10-2$. We consider the assignment $\mathcal{K}=\operatorname{diag}(S_0,I_0,I_0,I_0,I_0)$, with anchor $(i_0,j_0)=(1,4)$ and $K_T=\mu^{-1}$, which is the reference timescale for all death processes. Using \Cref{unimodularmain} and Algorithm C.4, we select a $\mathbb{Z}$-basis of $\Delta\Gamma$ that simplifies the coefficient exponents, giving the Laurent basis $\mathcal{G}=\{g_1,\dots,g_8\}$, where
\begin{equation}
\begin{aligned}
g_1&=\frac{\beta I_0}{\mu}, & g_2&=\frac{\Lambda}{\mu S_0}, & g_3&=\frac{\omega I_0}{\mu S_0}, & g_4&=\frac{\nu_V}{\mu}, & g_5&=\frac{\omega_V I_0}{\mu S_0}, & g_6&=\frac{\sigma}{\mu}, & g_7&=\frac{\gamma}{\mu}, & g_8&=\frac{S_0}{I_0}.
\end{aligned}
\end{equation}
The non-dimensionalised SEIRV scheme is thus expressible as
\begin{equation}
    \frac{d\mathbf{x}}{dt} =
    \underbrace{\begin{bmatrix}
-1 &  1 &  1 & -1 & -1 &  1 &  0 &  0 &  0 &  0 &  0 &  0 \\
 g_8 &  0 &  0 &  0 &  0 &  0 & -1 & -1 &  0 &  0 &  0 &  0 \\
 0 &  0 &  0 &  0 &  0 &  0 &  0 &  1 & -1 & -1 &  0 &  0 \\
 0 &  0 & -g_8 &  0 &  0 &  0 &  0 &  0 &  1 &  0 & -1 &  0 \\
 0 &  0 &  0 &  0 &  g_8 & -g_8 &  0 &  0 &  0 &  0 &  0 & -1
\end{bmatrix}}_{\tilde{S}}
\underbrace{\begin{bmatrix}
g_1 x_1 x_3 \\
g_2 \\
g_3 x_4 \\
x_1 \\
g_4 x_1 \\
g_5 x_5 \\
x_2 \\
g_6 x_2 \\
g_7 x_3 \\
x_3 \\
x_4 \\
x_5
\end{bmatrix}}_{\tilde{V}(\mathbf{x})},
\label{eq:seirvnondim}
\end{equation}
with $\mathbf{x}(0) = (1,0,1,0,0)^\intercal$. We have that for this configuration, $F_\mathrm{mono}=17$, $F_\mathrm{gen}=22$, and $F_\mathrm{time}=1$, and all monomial coefficients are polynomial in the generators so that $\mathcal{G}$ is also a polynomial basis. The first seven generators compare reaction fluxes at the reference populations with the corresponding death fluxes, while $g_8$ measures the initial susceptible-to-infectious ratio.
\subsection{Kim-Forger Model}\label{subsec:kimforger}
The Kim-Forger model arises from a simplified circadian clock gene regulatory network (GRN), characterised by the reactions
\begin{equation}
\begin{aligned}
\mathrm{S} &\xrightarrow{k_1} \mathrm{S} + \mathrm{X}, &
\mathrm{X} &\xrightarrow{k_2} \emptyset, &
\mathrm{X} &\xrightarrow{k_3} \mathrm{X} + \mathrm{Y}, &
\mathrm{Y} &\xrightarrow{k_4} \emptyset, & \mathrm{Y} &\xrightarrow{k_5} \mathrm{Y} + \mathrm{Z}, &
\mathrm{Z} &\xrightarrow{k_6} \emptyset, &
\mathrm{S} + \mathrm{Z} &\underset{k_b}{\stackrel{k_f}{\rightleftharpoons}} \mathrm{C}, &
\mathrm{C} &\xrightarrow{k_u} \mathrm{S},
\end{aligned} \label{eq:kimforgereqns}
\end{equation}
where $\mathrm{X}$, $\mathrm{Y}$, $\mathrm{Z}$, $\mathrm{S}$, and $\mathrm{C}$ denote the clock gene mRNA, encoded protein, repressor protein, activator/transcription factor, and sequestered activator-repressor complex, respectively \cite{Kim2012, Lapuz2026}. Denoting $\mathbf{X}=([\mathrm{X}],[\mathrm{Y}],[\mathrm{Z}],[\mathrm{S}],[\mathrm{C}])^\intercal \in \mathbb{R}^5_{\geq 0}$, we have that by the law of mass action, \eqref{eq:kimforgereqns} is expressible as
\begin{equation}
\begin{aligned}
\frac{dX_{1}}{dT} &= \fitRHS{k_1 X_{4} - k_2 X_{1}}, \\
\frac{dX_{2}}{dT} &= \fitRHS{k_3 X_{1} - k_4 X_{2}}, \\
\frac{dX_{3}}{dT} &= \fitRHS{k_5 X_{2} - k_6 X_{3} - k_f X_{3} X_{4} + k_b X_{5}}, \\
\frac{dX_{4}}{dT} &= \fitRHS{- k_f X_{3} X_{4} + k_b X_{5} + k_u X_{5}}, \\
\frac{dX_{5}}{dT} &= \fitRHS{k_f X_{3} X_{4} - k_b X_{5} - k_u X_{5}},
\end{aligned}
\implies
\begin{aligned}
\frac{d\mathbf{X}}{dT} &= \underbrace{\begin{bmatrix}
1 & -1 & 0 & 0 & 0 & 0 & 0 & 0 & 0 \\
0 & 0 & 1 & -1 & 0 & 0 & 0 & 0 & 0 \\
0 & 0 & 0 & 0 & 1 & -1 & -1 & 1 & 0 \\
0 & 0 & 0 & 0 & 0 & 0 & -1 & 1 & 1 \\
0 & 0 & 0 & 0 & 0 & 0 & 1 & -1 & -1
\end{bmatrix}}_{S}\underbrace{\begin{bmatrix}
k_1X_{4} \\
k_2X_{1} \\
k_3X_{1} \\
k_4X_{2} \\
k_5X_{2} \\
k_6X_{3} \\
k_fX_{3}X_{4} \\
k_bX_{5} \\
k_uX_{5}
\end{bmatrix}}_{V(\mathbf{X})}.
\end{aligned}
\label{eq:kimforgerdim}
\end{equation}
We consider initial conditions of $\mathbf{X}(0)=(x_0,y_0,z_0,s_0,0)^\intercal$, with typical numerical initial conditions in Table D.5. We see that there are nine distinct rate symbols and graph $G$ has components $\{\mathrm{X}\}$, $\{\mathrm{Y}\}$ and $\{\mathrm{Z},\mathrm{S},\mathrm{C}\}$, with initial-value sets $W_1=\{x_0\}$, $W_2=\{y_0\}$ and $W_3=\{z_0,s_0\}$. Consequently, $H_0$ has three components and \Cref{thm:ivpmin} gives $N^\ast_\mathrm{IVP}=9-1+4-3=9$, so that among all assignments from $U_\mathrm{cand}=\{x_0,y_0,z_0,s_0\}$, this minimum is attained exactly when $K_1=x_0$, $K_2=y_0$ and $K_3,K_4,K_5\in\{z_0,s_0\}$. Fixing $K_3=z_0$ and $K_4=s_0$, this reduces to $K_5\in\{z_0,s_0\}$, giving nine coefficient generators. Choosing $K_5\in\{x_0,y_0\}$ instead gives ten, while $K_1=s_0$ and $K_2=z_0$, with the other initial scalings retained, join all three components and give eleven IVP parameters. The Buckingham count $N_\Pi=13-2=11$ exceeds the minimum by two because of the independent rescalings of the $\mathrm{X}$ and $\mathrm{Y}$ components described in \Cref{prop:buckingham}.\\~\\
We execute Algorithm C.8 over the time window $[T_\mathrm{min}, T_\mathrm{max}]$, where $T_\mathrm{min}=6\mathrm{~hr}$ and $T_\mathrm{max}=12\mathrm{~hr}$, using the rate constants in Table D.4, $x_0=y_0=1\mathrm{~\mu M}$ and $z_0=s_0=0.06128\mathrm{~\mu M}$. The two remaining choices of $K_5$ are numerically indistinguishable under both deviation-based objectives, with $R^\ast \approx 1.967368$. Nonetheless, the optimal assignment was found to be $\sigma^{\mathrm{lex}}=(x_0,y_0,z_0,s_0,s_0)$, with anchor $(i_0,j_0)=(4,8)$ and $K_T=k_b^{-1}$, which is the reference timescale for $\mathrm{C}\xrightarrow{k_b}\mathrm{S}+\mathrm{Z}$. Using \Cref{unimodularmain} and Algorithm C.4, we select a $\mathbb{Z}$-basis of $\Delta\Gamma$ that simplifies the coefficient exponents, giving the Laurent basis $\mathcal{G}=\{g_1,\dots,g_9\}$, where
\begin{equation}
\begin{aligned}
g_1 &= \frac{k_1 s_0}{k_b x_0}, &
g_2 &= \frac{k_2}{k_b}, &
g_3 &= \frac{k_3 x_0}{k_b y_0}, &
g_4 &= \frac{k_4}{k_b}, &
g_5 &= \frac{k_5 y_0}{k_b z_0}, & g_6 &= \frac{k_6}{k_b}, &
g_7 &= \frac{k_f z_0}{k_b}, &
g_8 &= \frac{k_u}{k_b}, &
g_9 &= \frac{s_0}{z_0},
\end{aligned}
\end{equation}
so that the non-dimensionalised Kim-Forger scheme is expressible as 
\begin{equation}
\resizeboxalt{0.94\linewidth}{!}{$\displaystyle
\frac{d\mathbf{x}}{dt} =
    \underbrace{\begin{bmatrix}
 1 & -1 &  0 &  0 &  0 &  0 &   0 &   0 &   0 \\
 0 &  0 &  1 & -1 &  0 &  0 &   0 &   0 &   0 \\
 0 &  0 &  0 &  0 &  1 & -1 & -g_9 &  g_9 &   0 \\
 0 &  0 &  0 &  0 &  0 &  0 &  -1 &   1 &   1 \\
 0 &  0 &  0 &  0 &  0 &  0 &   1 &  -1 &  -1
\end{bmatrix}}_{\tilde{S}}
\underbrace{\begin{bmatrix}
g_1 x_4 \\
g_2 x_1 \\
g_3 x_1 \\
g_4 x_2 \\
g_5 x_2 \\
g_6 x_3 \\
g_7 x_3 x_4 \\
x_5 \\
g_8 x_5
\end{bmatrix}}_{\tilde{V}(\mathbf{x})},
$}
\label{eq:kimforgernondim}
\end{equation}
with $\mathbf{x}(0) = (1,1,1,1,0)^\intercal$. Clearly, the non-dimensionalised system has $9$ free parameters, with $F_\mathrm{mono}=13$, $F_\mathrm{gen}=25$, and $F_\mathrm{time}=1$, and all monomial coefficients are polynomial in the generators so that $\mathcal{G}$ is also a polynomial basis. The production parameters $g_1$, $g_3$ and $g_5$ incorporate the chosen species scales, while $g_2$, $g_4$, $g_6$ and $g_8$ compare first-order losses with unbinding, while the association parameter $g_7$ compares the characteristic sequestration frequency $k_fz_0$ with $k_b$, and $g_9$ measures the initial activator-to-repressor ratio.
\subsection{Vilar's Model of Genetic Oscillators}\label{subsec:vilar}
Vilar's model is a minimal model of genomically based genetic oscillation that describes the interactions between an activator and a repressor gene, capturing noise-resistant circadian-like oscillations in transcriptional feedback loops \cite{Vilar2002}. The model is characterised by the reactions
\pagebreak[0]
\begin{equation}
\resizeboxalt{0.9\linewidth}{!}{$\displaystyle
\begin{aligned}
\mathrm{DA} &\xrightarrow{\alpha_A} \mathrm{DA} + \mathrm{MA}, &
\mathrm{DA^\prime} &\xrightarrow{\alpha_{A^\prime}} \mathrm{DA^\prime} + \mathrm{MA}, &
\mathrm{DR} &\xrightarrow{\alpha_R} \mathrm{DR} + \mathrm{MR}, &
\mathrm{DR^\prime} &\xrightarrow{\alpha_{R^\prime}} \mathrm{DR^\prime} + \mathrm{MR}, \\
\mathrm{MA} &\xrightarrow{\beta_A} \mathrm{A} + \mathrm{MA}, &
\mathrm{MR} &\xrightarrow{\beta_R} \mathrm{MR} + \mathrm{R}, &
\mathrm{A} + \mathrm{DA} &\xrightarrow{\gamma_A} \mathrm{DA^\prime}, &
\mathrm{A} + \mathrm{R} &\xrightarrow{\gamma_C} \mathrm{C}, \\
\mathrm{A} + \mathrm{DR} &\xrightarrow{\gamma_R} \mathrm{DR^\prime}, &
\mathrm{A} &\xrightarrow{\delta_A} \emptyset, &
\mathrm{C} &\xrightarrow{\delta_A} \mathrm{R}, &
\mathrm{MA} &\xrightarrow{\delta_{MA}} \emptyset, \\
\mathrm{MR} &\xrightarrow{\delta_{MR}} \emptyset, &
\mathrm{R} &\xrightarrow{\delta_R} \emptyset, &
\mathrm{DA^\prime} &\xrightarrow{\theta_A} \mathrm{A} + \mathrm{DA}, &
\mathrm{DR^\prime} &\xrightarrow{\theta_R} \mathrm{A} + \mathrm{DR},
\end{aligned}$}
\label{eq:vilar}
\end{equation}
where species are as in \Cref{tab:vilar2002_species}.
\begin{table}[H]
\centering
\resizeboxalt{\textwidth}{!}{$
\begin{tabular}{|l|l|l|l|l|l|}\hline
\textbf{Name} & \textbf{Description} & \textbf{Name} & \textbf{Description} & \textbf{Name} & \textbf{Description} \\ \hline
$\mathrm{A}$ & activator protein & $\mathrm{C}$ & A:R inactive complex & $\mathrm{DA}$ & unbound activator gene promoter \\
$\mathrm{DA^\prime}$ & activator promoter bound by A & $\mathrm{DR}$ & unbound repressor gene promoter & $\mathrm{DR^\prime}$ & repressor promoter bound by A \\
$\mathrm{MA}$ & activator mRNA & $\mathrm{MR}$ & repressor mRNA & $\mathrm{R}$ & repressor protein \\
\hline\end{tabular}$}
\caption{Involved species in the model.}
\label{tab:vilar2002_species}
\end{table}
Denoting $\mathbf{X}:=([\mathrm{A}],[\mathrm{C}],[\mathrm{DA}],[\mathrm{DA^\prime}],[\mathrm{DR}],[\mathrm{DR^\prime}],[\mathrm{MA}],[\mathrm{MR}],[\mathrm{R}])^\intercal \in \mathbb{R}^9_{\geq 0}$, we have, by the law of mass action, that the corresponding dimensional system for \eqref{eq:vilar} has the factorisation
\begin{equation}
\frac{d\mathbf{X}}{dT} = \underbrace{\begin{bmatrix}
0 & 0 & 0 & 0 & 1 & 0 & -1 & -1 & -1 & -1 & 0 & 0 & 0 & 0 & 1 & 1 \\
0 & 0 & 0 & 0 & 0 & 0 & 0 & 1 & 0 & 0 & -1 & 0 & 0 & 0 & 0 & 0 \\
0 & 0 & 0 & 0 & 0 & 0 & -1 & 0 & 0 & 0 & 0 & 0 & 0 & 0 & 1 & 0 \\
0 & 0 & 0 & 0 & 0 & 0 & 1 & 0 & 0 & 0 & 0 & 0 & 0 & 0 & -1 & 0 \\
0 & 0 & 0 & 0 & 0 & 0 & 0 & 0 & -1 & 0 & 0 & 0 & 0 & 0 & 0 & 1 \\
0 & 0 & 0 & 0 & 0 & 0 & 0 & 0 & 1 & 0 & 0 & 0 & 0 & 0 & 0 & -1 \\
1 & 1 & 0 & 0 & 0 & 0 & 0 & 0 & 0 & 0 & 0 & -1 & 0 & 0 & 0 & 0 \\
0 & 0 & 1 & 1 & 0 & 0 & 0 & 0 & 0 & 0 & 0 & 0 & -1 & 0 & 0 & 0 \\
0 & 0 & 0 & 0 & 0 & 1 & 0 & -1 & 0 & 0 & 1 & 0 & 0 & -1 & 0 & 0
\end{bmatrix}}_{S}
\underbrace{\begin{bmatrix}
\alpha_A X_{3} \\
\alpha_{A^\prime} X_{4} \\
\alpha_R X_{5} \\
\alpha_{R^\prime} X_{6} \\
\beta_A X_{7} \\
\beta_R X_{8} \\
\gamma_A X_{1}X_{3} \\
\gamma_C X_{1}X_{9} \\
\gamma_R X_{1}X_{5} \\
\delta_A X_{1} \\
\delta_A X_{2} \\
\delta_{MA} X_{7} \\
\delta_{MR} X_{8} \\
\delta_R X_{9} \\
\theta_A X_{4} \\
\theta_R X_{6}
\end{bmatrix}}_{V(\mathbf{X})}. \label{eq:vilar2002dim}
\end{equation}
This has fifteen distinct rate constants, with typical numerical values given in Table D.7, and we consider initial conditions of $X_3(0) = [\mathrm{DA}]_0$ and $X_5(0) = [\mathrm{DR}]_0$ and all other species initially zero, with typical numerical initial conditions in Table D.8. The graph $G$ has components $\{1,2,3,4,5,6,9\}$, $\{7\}$ and $\{8\}$. Both positive initial values lie in the first component, so that $H_0$ consists of the edge joining $[\mathrm{DA}]_0$ and $[\mathrm{DR}]_0$, and the shared rate $\delta_A$ occurs in $\mathrm{A}\to\emptyset$ and $\mathrm{C}\to\mathrm{R}$, whose reactants and active species lie in the same component of $G$, so that their vectors $\bar{z}_j$ therefore vanish. Thus, every assignment from $U_\mathrm{cand}=\{[\mathrm{DA}]_0,[\mathrm{DR}]_0\}$ attains $N^\ast_\mathrm{IVP}=15-1+2-1=15$, equal to the Buckingham count $17-2$.\\~\\
Accordingly fixing $K_3=[\mathrm{DA}]_0$ and $K_5=[\mathrm{DR}]_0$, and setting $U_\mathrm{cand}:=\set{[\mathrm{DA}]_0,[\mathrm{DR}]_0}$, we execute Algorithm C.8 over the time window $[T_\mathrm{min}, T_\mathrm{max}]$, where $T_\mathrm{min}=100\mathrm{~hr}$ and $T_\mathrm{max}=400\mathrm{~hr}$. Since $[\mathrm{DA}]_0=[\mathrm{DR}]_0$, all generator-minimising assignments are indistinguishable under both the supremum and lexicographic deviation-based objectives, with $R^\ast\approx7.688287$. Nonetheless, the optimal assignment was found to be $\sigma^{\mathrm{lex}}=([\mathrm{DR}]_0,\allowbreak[\mathrm{DR}]_0,\allowbreak[\mathrm{DA}]_0,\allowbreak[\mathrm{DA}]_0,\allowbreak[\mathrm{DR}]_0,\allowbreak[\mathrm{DR}]_0,\allowbreak[\mathrm{DA}]_0,\allowbreak[\mathrm{DR}]_0,\allowbreak[\mathrm{DR}]_0)$, with anchor $(i_0,j_0)=(1,10)$ and $K_T=\delta_A^{-1}$, corresponding to the reference timescale for $\mathrm{A}\xrightarrow{\delta_A}\emptyset$. Using \Cref{unimodularmain} and Algorithm C.4, we select a $\mathbb{Z}$-basis of $\Delta\Gamma$ that simplifies the coefficient exponents, giving the Laurent basis $\mathcal{G}=\{g_1,\dots,g_{15}\}$, where
\begin{equation}
\begin{aligned}
  g_1  &= \frac{\alpha_A}{\delta_A}, &
  g_2  &= \frac{\alpha_{A^\prime}}{\delta_A}, &
  g_3  &= \frac{\alpha_R}{\delta_A}, &
  g_4  &= \frac{\alpha_{R^\prime}}{\delta_A}, &
  g_5  &= \frac{\beta_A[\mathrm{DA}]_0}{\delta_A[\mathrm{DR}]_0}, \\
  g_6  &= \frac{\beta_R}{\delta_A}, &
  g_7  &= \frac{\gamma_A[\mathrm{DR}]_0}{\delta_A}, &
  g_8  &= \frac{\gamma_C[\mathrm{DR}]_0}{\delta_A}, &
  g_9  &= \frac{\gamma_R[\mathrm{DR}]_0}{\delta_A}, &
  g_{10} &= \frac{\delta_{MA}}{\delta_A}, \\
  g_{11} &= \frac{\delta_{MR}}{\delta_A}, &
  g_{12} &= \frac{\delta_R}{\delta_A}, &
  g_{13} &= \frac{\theta_A}{\delta_A}, &
  g_{14} &= \frac{\theta_R}{\delta_A}, &
  g_{15} &= \frac{[\mathrm{DA}]_0}{[\mathrm{DR}]_0},
    \end{aligned}
    \end{equation}
so that the non-dimensionalised model is expressible as 
\begin{equation}
\resizeboxalt{0.94\linewidth}{!}{$\displaystyle
\frac{d\mathbf{x}}{dt} =
    \underbrace{\begin{bmatrix}
      0 & 0 & 0 & 0 & 1 & 0 & -g_{15} & -1 & -1 & -1 &  0 &  0 &  0 &  0 &  g_{15} &  1 \\
      0 & 0 & 0 & 0 & 0 & 0 &  0 &  1 &  0 &  0 & -1 &  0 &  0 &  0 &  0 &  0 \\
      0 & 0 & 0 & 0 & 0 & 0 & -1 &  0 &  0 &  0 &  0 &  0 &  0 &  0 &  1 &  0 \\
      0 & 0 & 0 & 0 & 0 & 0 &  1 &  0 &  0 &  0 &  0 &  0 &  0 &  0 & -1 &  0 \\
      0 & 0 & 0 & 0 & 0 & 0 &  0 &  0 & -1 &  0 &  0 &  0 &  0 &  0 &  0 &  1 \\
      0 & 0 & 0 & 0 & 0 & 0 &  0 &  0 &  1 &  0 &  0 &  0 &  0 &  0 &  0 & -1 \\
      1 & 1 & 0 & 0 & 0 & 0 &  0 &  0 &  0 &  0 &  0 & -1 &  0 &  0 &  0 &  0 \\
      0 & 0 & 1 & 1 & 0 & 0 &  0 &  0 &  0 &  0 &  0 &  0 & -1 &  0 &  0 &  0 \\
      0 & 0 & 0 & 0 & 0 & 1 &  0 & -1 &  0 &  0 &  1 &  0 &  0 & -1 &  0 &  0
\end{bmatrix}}_{\tilde{S}}
\underbrace{\begin{bmatrix}
    g_1 x_3 \\
    g_2 x_4 \\
    g_3 x_5 \\
    g_4 x_6 \\
    g_5 x_7 \\
    g_6 x_8 \\
    g_7 x_1 x_3 \\
    g_8 x_1 x_9 \\
    g_9 x_1 x_5 \\
    x_1 \\
    x_2 \\
    g_{10} x_7 \\
    g_{11} x_8 \\
    g_{12} x_9 \\
    g_{13} x_4 \\
    g_{14} x_6
    \end{bmatrix}}_{\tilde{V}(\mathbf{x})},
$}
\label{eq:vilarnondim}
\end{equation}
with $\mathbf{x}(0) = (0,0,1,0,1,0,0,0,0)^\intercal$. Clearly, the non-dimensionalised system has $15$ free parameters, with $F_\mathrm{mono}=26$, $F_\mathrm{gen}=35$, and $F_\mathrm{time}=1$, and all monomial coefficients are polynomial in the generators so that $\mathcal{G}$ is also a polynomial basis. The generators compare transcription, translation, association, degradation and promoter dissociation with the activator degradation rate, including the prescribed species scales where required, while $g_{15}$ measures the initial activator-to-repressor promoter ratio.\\~\\
The structural quantities and representation scores for each example are summarised in \Cref{tab:summary}. As explained in \Cref{rem:scope}, allowing rate-dependent monomial scales for the initially absent mRNAs in Vilar's model gives
\begin{equation}\label{eq:vilar-monomial-scales}
K_\mathrm{MA}=\frac{\alpha_A[\mathrm{DA}]_0}{\delta_A}, K_\mathrm{MR}=\frac{\alpha_R[\mathrm{DR}]_0}{\delta_A},
(g_1,g_2,g_5)\mapsto
\left(\frac{\alpha_{A^\prime}}{\alpha_A},\frac{\beta_A\alpha_A[\mathrm{DA}]_0}{\delta_A^2[\mathrm{DR}]_0}\right), 
(g_3,g_4,g_6)\mapsto
\left(\frac{\alpha_{R^\prime}}{\alpha_R},\frac{\beta_R\alpha_R}{\delta_A^2}\right).
\end{equation}
With all other scales unchanged, the remaining nine parameters are unaltered, yielding thirteen algebraically independent combinations in which every coefficient is polynomial. More details are available in Section D.5.1, however, a characterisation of parameter-minimal assignments within this larger class is beyond the present results.
\begin{table}[t]
\centering
\resizeboxalt{\textwidth}{!}{$
\begin{tabular}{|l|c|c|c|c|c|c|c|c|c|}\hline
\textbf{Model} & $n$ & $r$ & $h$ & $q$ & $|\iota(\mathcal{J})|$ & $c(H_0)$ & $N^\ast_\mathrm{IVP}$ & $N_\Pi$ & $(F_\mathrm{mono},F_\mathrm{gen},F_\mathrm{time})$ \\ \hline
Michaelis--Menten, first assignment & 4 & 4 & 4 & 1 & 2 & 1 & 4 & 4 & (11,10,1) \\
Michaelis--Menten, second assignment & 4 & 4 & 4 & 1 & 2 & 1 & 4 & 4 & (13,10,1) \\
SEIRV & 5 & 12 & 8 & 1 & 2 & 1 & 8 & 8 & (17,22,1) \\
Kim--Forger & 5 & 9 & 9 & 3 & 4 & 3 & 9 & 11 & (13,25,1) \\
Vilar & 9 & 16 & 15 & 3 & 2 & 1 & 15 & 15 & (26,35,1) \\
\hline\end{tabular}$}
\caption{Structural quantities and scores for the reported non-dimensionalisations. The numbers of species, reactions, distinct rate symbols and components of $G$ are denoted by $n$, $r$, $h$ and $q$, respectively. The initial-value graph $H_0$ determines $N^\ast_\mathrm{IVP}$ through \Cref{thm:ivpmin}, while $N_\Pi$ is the Buckingham count of \Cref{prop:buckingham}.}
\label{tab:summary}
\end{table}

\section{Conclusion}
In this work, we present a unified algebraic framework that transforms the non-dimensionalisation of polynomial dynamical systems from a largely intuition-based procedure into a rigorous structural problem by viewing them as CRNs. By mapping admissible scalings to a subgroup of the corresponding Laurent polynomial ring, we established that the theoretical minimum number of multiplicative dimensionless parameters for a fixed reaction-channel coefficient set and scale assignment, with an appropriate time scale, is $h+s-c_\mathrm{scale}+r_\mathrm{rate}-1$. On the other hand, including the nonzero dimensionless initial conditions gives the full-IVP count $h+s_\mathrm{tot}-c_\mathrm{aug}+r_\mathrm{IVP}-1$, with these formulae explicitly linking parameter dependence to the connectivity of the associated scale graphs and the stoichiometric rate lattice. In particular, scaling positive initial concentrations to unity can attain the structural lower bound, whereas a smaller coefficient count may merely transfer parameters into the initial conditions. The IVP count does not exceed the corresponding Buckingham count, with any strict reduction reflecting scaling freedom permitted by the equations beyond a common change of physical units. This result also provides a precise answer to the classical modelling problem of determining exactly when characteristic scales can be introduced without increasing the number of free parameters. Furthermore, we established that, for a given non-dimensionalisation, distinct bases of the same coefficient lattice are related by unimodular transformations, revealing that they are just alternative coordinatisations of the same intrinsic exponent lattice.\\~\\
A distinctive feature of our approach is that it decouples two conceptually different aspects of non-dimensionalisation, namely the structural minimality of the free parameter set and the numerical scaling of the variables. As demonstrated by our implementation in CRNNondim, this separation allows for a three-stage workflow in which we first determine the space of parameter-minimal admissible assignments, followed by employing discrete optimisation to select those minimising the logarithmic trajectory deviations, and lastly, unimodular simplification to select an optimal non-dimensionalisation for a given CRN, and thus for the corresponding polynomial system. Finally, we established the computational efficiency of the framework through a rigorous analysis of the algorithms employed by CRNNondim, ensuring its applicability to high-dimensional systems pervasive across the quantitative life and physical sciences.\\~\\
It is therefore natural to extend this framework beyond polynomial dynamical systems to the broader class of rational dynamical systems \cite{cornish1979}, as well as to systems with spatial dependence \cite{Murray2003,Turing1952}, such as partial differential equations. Allowing species scales to depend monomially on the rate constants also enlarges the admissible class and can lower the parameter count, as demonstrated for the zero-initial-value mRNA species in Vilar's model. Another avenue for future research is the refinement of the unimodular simplification algorithms, where local greedy searches currently lead to sparse bases in practice, but the underlying optimisation over $\mathrm{GL}_{d-1}(\mathbb{Z})$ suggests more sophisticated heuristics could be developed using tools from lattice basis reduction \cite{Lenstra1982}. Ultimately, by making non-dimensionalisation automatic, rigorous, and systematic, CRNNondim is positioned to become integrated into the routine analysis of polynomial dynamical systems, further solidifying the role of CRNs as a practical framework for representing and analysing these complex models.
\par
\begin{samepage}
\section{CRediT authorship contribution statement}
\textbf{Georgio Hawi}: conceptualisation, data curation, formal analysis, investigation, methodology, project administration, resources, software, validation, visualisation, writing --- original draft, writing --- review \& editing. \\
\textbf{Martin Wechselberger}: conceptualisation, formal analysis, funding acquisition, investigation, methodology, project administration, resources, supervision, validation, visualisation, writing --- original draft, writing --- review \& editing. 
\par
\end{samepage}
\section{Declaration of Competing Interests}
The authors declare that they have no known competing financial interests or personal relationships that could have appeared to influence the work reported in this paper.
\section{Data availability}
All data and procedures are available within the manuscript and its Supporting Information file. Source code and usage instructions for CRNNondim are available at \url{https://github.com/georgio-hawi/crnnondim}.
\section{Acknowledgements}
GH and MW gratefully acknowledge support from the Australian Research Council Discovery Projects (DP220101817 and DP260100522).
\putbib[References]
\end{bibunit}


\begin{thebibliography}{10}

\bibitem{Strogatz2018}
Strogatz SH.
\newblock {N}onlinear {D}ynamics and {C}haos.
\newblock CRC Press; 2018.
\newblock Available from: \url{https://doi.org/10.1201/9780429492563}.

\bibitem{Hethcote2000}
Hethcote HW.
\newblock {T}he {M}athematics of {I}nfectious {D}iseases.
\newblock SIAM Review. 2000 Dec;42(4):599-653.
\newblock Available from: \url{https://doi.org/10.1137/S0036144500371907}.

\bibitem{Lorenz2004}
Lorenz EN.
\newblock {D}eterministic {N}onperiodic {F}low.
\newblock In: The Theory of Chaotic Attractors. Springer New York; 2004. p. 25-36.
\newblock Available from: \url{https://doi.org/10.1007/978-0-387-21830-4_2}.

\bibitem{Volterra1926}
Volterra V.
\newblock {F}luctuations in the {A}bundance of a {S}pecies considered {M}athematically.
\newblock Nature. 1926 Oct;118(2972):558-60.
\newblock Available from: \url{https://doi.org/10.1038/118558a0}.

\bibitem{Horn1972}
Horn F, Jackson R.
\newblock {G}eneral mass action kinetics.
\newblock Archive for Rational Mechanics and Analysis. 1972 Jan;47(2):81-116.
\newblock Available from: \url{https://doi.org/10.1007/BF00251225}.

\bibitem{kermack1927contribution}
Kermack WO, McKendrick AG.
\newblock {A} contribution to the mathematical theory of epidemics.
\newblock Proceedings of the Royal Society of London Series A, Containing papers of a mathematical and physical character. 1927;115(772):700-21.

\bibitem{Nagumo1962}
Nagumo J, Arimoto S, Yoshizawa S.
\newblock {A}n {A}ctive {P}ulse {T}ransmission {L}ine {S}imulating {N}erve {A}xon.
\newblock Proceedings of the IRE. 1962 Oct;50(10):2061-70.
\newblock Available from: \url{https://doi.org/10.1109/JRPROC.1962.288235}.

\bibitem{Buckingham1914}
Buckingham E.
\newblock {O}n {P}hysically {S}imilar {S}ystems; {I}llustrations of the {U}se of {D}imensional {E}quations.
\newblock Physical Review. 1914 Oct;4(4):345-76.
\newblock Available from: \url{https://doi.org/10.1103/PhysRev.4.345}.

\bibitem{Barenblatt2003}
Barenblatt GI.
\newblock {S}caling.
\newblock Cambridge University Press; 2003.
\newblock Available from: \url{https://doi.org/10.1017/CBO9780511814921}.

\bibitem{Curtis1982}
Curtis WD, Logan JD, Parker WA.
\newblock {D}imensional analysis and the pi theorem.
\newblock Linear Algebra and its Applications. 1982 Oct;47:117-26.
\newblock Available from: \url{https://doi.org/10.1016/0024-3795(82)90229-4}.

\bibitem{Feliu2022}
Feliu E, Walcher S, Wiuf C.
\newblock {C}ritical {P}arameters for {S}ingular {P}erturbation {R}eductions of {C}hemical {R}eaction {N}etworks.
\newblock Journal of Nonlinear Science. 2022 Sep;32(6):83.
\newblock Available from: \url{https://doi.org/10.1007/s00332-022-09843-4}.

\bibitem{Goeke2015}
Goeke A, Walcher S, Zerz E.
\newblock {D}etermining ``small parameters'' for quasi-steady state.
\newblock Journal of Differential Equations. 2015 Aug;259(3):1149-80.
\newblock Available from: \url{https://doi.org/10.1016/j.jde.2015.02.038}.

\bibitem{Goeke2013}
Goeke A, Walcher S.
\newblock {Q}uasi-{S}teady {S}tate: {S}earching for and {U}tilizing {S}mall {P}arameters.
\newblock In: Recent Trends in Dynamical Systems. Springer Basel; 2013. p. 153-78.
\newblock Available from: \url{https://doi.org/10.1007/978-3-0348-0451-6_8}.

\bibitem{Lapuz2026}
Lapuz TEF, Wechselberger M.
\newblock {C}oordinate-{I}ndependent {M}odel {R}eductions of {C}hemical {R}eaction {N}etworks {B}ased on {G}eometric {S}ingular {P}erturbation {T}heory.
\newblock SIAM Journal on Life Sciences. 2026;1(2):163-201.
\newblock Available from: \url{https://epubs.siam.org/doi/abs/10.1137/25M1785022}.

\bibitem{Hairer1996}
Hairer E, Wanner G.
\newblock {S}olving {O}rdinary {D}ifferential {E}quations {II}.
\newblock Springer Berlin Heidelberg; 1996.
\newblock Available from: \url{https://doi.org/10.1007/978-3-642-05221-7}.

\bibitem{bridgman1922dimensional}
Bridgman PW.
\newblock {D}imensional {A}nalysis.
\newblock New Haven : Yale University Press; 1922.

\bibitem{murray2002mathematical}
Murray JD.
\newblock {M}athematical {B}iology: {I}. {A}n {I}ntroduction.
\newblock Springer New York; 2002.
\newblock Available from: \url{https://doi.org/10.1007/b98868}.

\bibitem{Lin1988}
Lin CC, Segel LA.
\newblock {M}athematics {A}pplied to {D}eterministic {P}roblems in the {N}atural {S}ciences.
\newblock Society for Industrial and Applied Mathematics; 1988.
\newblock Available from: \url{https://doi.org/10.1137/1.9781611971347}.

\bibitem{szirtes2007applied}
Szirtes T.
\newblock {A}pplied {D}imensional {A}nalysis and {M}odeling.
\newblock Butterworth-Heinemann; 2007.

\bibitem{Hart1995}
Hart GW.
\newblock {M}ultidimensional {A}nalysis.
\newblock Springer New York; 1995.
\newblock Available from: \url{https://doi.org/10.1007/978-1-4612-4208-6}.

\bibitem{HubertSedoglavic2006}
Hubert {\'E}, Sedoglavic A.
\newblock {P}olynomial {T}ime {N}ondimensionalisation of {O}rdinary {D}ifferential {E}quations via {T}heir {L}ie {P}oint {S}ymmetries; 2006.
\newblock Preprint, arXiv:cs/0604060.
\newblock Available from: \url{https://arxiv.org/abs/cs/0604060}.

\bibitem{HubertLabahn2013}
Hubert {\'E}, Labahn G.
\newblock {S}caling {I}nvariants and {S}ymmetry {R}eduction of {D}ynamical {S}ystems.
\newblock Foundations of Computational Mathematics. 2013;13(4):479-516.

\bibitem{Atherton2014}
Atherton MA, Bates RA, Wynn HP.
\newblock {D}imensional {A}nalysis {U}sing {T}oric {I}deals: {P}rimitive {I}nvariants.
\newblock PLOS ONE. 2014;9(12):e112827.

\bibitem{Tanburn2025}
Tanburn R, Hendron D, Maini P, Amethyst S, Dufresne E, Harrington HA.
\newblock {N}ondimensionalization {I}s {M}ore {S}cience {T}han {A}rt; 2025.
\newblock Preprint, arXiv:2512.13455.
\newblock Available from: \url{https://arxiv.org/abs/2512.13455}.

\bibitem{hars1981inverse}
H{\'a}rs V, T{\'o}th J.
\newblock {O}n the inverse problem of reaction kinetics.
\newblock In: Qualitative Theory of Differential Equations. vol.~30; 1981. p. 363-79.

\bibitem{Craciun2020}
Craciun G, Johnston MD, Szederk{\'e}nyi G, Tonello E, Tóth J, Yu PY.
\newblock {R}ealizations of kinetic differential equations.
\newblock Mathematical Biosciences and Engineering. 2020;17(1):862-92.
\newblock Available from: \url{https://doi.org/10.3934/mbe.2020046}.

\bibitem{Feinberg1987}
Feinberg M.
\newblock {C}hemical reaction network structure and the stability of complex isothermal reactors---{I}. The deficiency zero and deficiency one theorems.
\newblock Chemical Engineering Science. 1987;42(10):2229-68.
\newblock Available from: \url{https://doi.org/10.1016/0009-2509(87)80099-4}.

\bibitem{Feinberg1988}
Feinberg M.
\newblock {C}hemical reaction network structure and the stability of complex isothermal reactors—{II}. {M}ultiple steady states for networks of deficiency one.
\newblock Chemical Engineering Science. 1988;43(1):1-25.
\newblock Available from: \url{https://doi.org/10.1016/0009-2509(88)87122-7}.

\bibitem{Feinberg2019}
Feinberg M.
\newblock {F}oundations of {C}hemical {R}eaction {N}etwork {T}heory.
\newblock Springer International Publishing; 2019.
\newblock Available from: \url{https://doi.org/10.1007/978-3-030-03858-8}.

\bibitem{Craciun2013}
Craciun G, Nazarov F, Pantea C.
\newblock {P}ersistence and {P}ermanence of {M}ass-{A}ction and {P}ower-{L}aw {D}ynamical {S}ystems.
\newblock SIAM Journal on Applied Mathematics. 2013 Jan;73(1):305-29.
\newblock Available from: \url{https://doi.org/10.1137/100812355}.

\bibitem{Godsil2001}
Godsil C, Royle G.
\newblock {A}lgebraic {G}raph {T}heory.
\newblock Springer New York; 2001.
\newblock Available from: \url{https://doi.org/10.1007/978-1-4613-0163-9}.

\bibitem{Biggs1974}
Biggs N.
\newblock {A}lgebraic {G}raph {T}heory.
\newblock Cambridge University Press; 1974.

\bibitem{Lang2002}
Lang S.
\newblock {A}lgebra.
\newblock Springer New York; 2002.
\newblock Available from: \url{https://doi.org/10.1007/978-1-4613-0041-0}.

\bibitem{Wegscheider1901}
Wegscheider R.
\newblock {\"U}ber simultane Gleichgewichte und die Beziehungen zwischen Thermodynamik und Reactionskinetik homogener Systeme.
\newblock Monatshefte f{\"u}r Chemie und verwandte Teile anderer Wissenschaften. 1901;22(8):849-906.
\newblock Available from: \url{https://doi.org/10.1007/BF01517498}.

\bibitem{osti_5633555}
McCormick ST.
\newblock {C}ombinatorial approach to some sparse-matrix problems.
\newblock Stanford Univ., CA (USA). Dept. of Operations Research; 1983.
\newblock Available from: \url{https://www.osti.gov/biblio/5633555}.

\bibitem{Meurer2017}
Meurer A, Smith CP, Paprocki M, {\v C}ert{\'i}k O, Kirpichev SB, Rocklin M, et~al.
\newblock {S}ymPy: symbolic computing in {P}ython.
\newblock PeerJ Computer Science. 2017 Jan;3:e103.
\newblock Available from: \url{https://doi.org/10.7717/peerj-cs.103}.

\bibitem{Harris2020}
Harris CR, Millman KJ, van~der Walt SJ, Gommers R, Virtanen P, Cournapeau D, et~al.
\newblock {A}rray programming with {NumPy}.
\newblock Nature. 2020 Sep;585(7825):357-62.
\newblock Available from: \url{https://doi.org/10.1038/s41586-020-2649-2}.

\bibitem{Virtanen2020}
Virtanen P, Gommers R, Oliphant TE, Haberland M, Reddy T, Cournapeau D, et~al.
\newblock {SciPy} 1.0: fundamental algorithms for scientific computing in {P}ython.
\newblock Nature Methods. 2020 Feb;17(3):261-72.
\newblock Available from: \url{https://doi.org/10.1038/s41592-019-0686-2}.

\bibitem{Hairer1999}
Hairer E, Wanner G.
\newblock {S}tiff differential equations solved by {R}adau methods.
\newblock Journal of Computational and Applied Mathematics. 1999 Nov;111(1--2):93-111.
\newblock Available from: \url{https://doi.org/10.1016/S0377-0427(99)00134-X}.

\bibitem{michaelis1913kinetik}
Michaelis L, Menten ML.
\newblock {D}ie {K}inetik der {I}nvertinwirkung.
\newblock Biochem Z. 1913;49:333-69.

\bibitem{Johnson2011}
Johnson KA, Goody RS.
\newblock {T}he {O}riginal {M}ichaelis {C}onstant: {T}ranslation of the 1913 {M}ichaelis--{M}enten {P}aper.
\newblock Biochemistry. 2011 Sep;50(39):8264-9.
\newblock Available from: \url{https://doi.org/10.1021/bi201284u}.

\bibitem{CornishBowden2015}
Cornish-Bowden A.
\newblock {O}ne hundred years of {M}ichaelis--{M}enten kinetics.
\newblock Perspectives in Science. 2015 Mar;4:3-9.
\newblock Available from: \url{https://doi.org/10.1016/j.pisc.2014.12.002}.

\bibitem{Anderson1991}
Anderson RM, May RM.
\newblock {I}nfectious {D}iseases of {H}umans: {D}ynamics and {C}ontrol.
\newblock Oxford University Press; 1991.
\newblock Available from: \url{https://doi.org/10.1093/oso/9780198545996.001.0001}.

\bibitem{Poonia2022}
Poonia RC, Saudagar AKJ, Altameem A, Alkhathami M, Khan MB, Hasanat MHA.
\newblock {A}n {E}nhanced {SEIR} {M}odel for {P}rediction of {COVID}-19 with {V}accination {E}ffect.
\newblock Life. 2022 Apr;12(5):647.
\newblock Available from: \url{https://doi.org/10.3390/life12050647}.

\bibitem{Kim2012}
Kim JK, Forger DB.
\newblock {A} mechanism for robust circadian timekeeping via stoichiometric balance.
\newblock Molecular Systems Biology. 2012 Dec;8(1):630.
\newblock Available from: \url{https://doi.org/10.1038/msb.2012.62}.

\bibitem{Vilar2002}
Vilar JMG, Kueh HY, Barkai N, Leibler S.
\newblock {M}echanisms of noise-resistance in genetic oscillators.
\newblock Proceedings of the National Academy of Sciences. 2002 Apr;99(9):5988-92.
\newblock Available from: \url{https://doi.org/10.1073/pnas.092133899}.

\bibitem{cornish1979}
Cornish-Bowden A.
\newblock {F}undamentals of {E}nzyme {K}inetics.
\newblock Butterworth-Heinemann; 1979.
\newblock Available from: \url{https://doi.org/10.1016/C2013-0-04130-8}.

\bibitem{Murray2003}
Murray JD.
\newblock {M}athematical {B}iology: {II}: {S}patial {M}odels and {B}iomedical {A}pplications.
\newblock Springer New York; 2003.
\newblock Available from: \url{https://doi.org/10.1007/b98869}.

\bibitem{Turing1952}
Turing AM.
\newblock {T}he chemical basis of morphogenesis.
\newblock Philosophical Transactions of the Royal Society of London Series B, Biological Sciences. 1952 Aug;237(641):37-72.
\newblock Available from: \url{https://doi.org/10.1098/rstb.1952.0012}.

\bibitem{Lenstra1982}
Lenstra AK, Lenstra HW, Lov{\'a}sz L.
\newblock {F}actoring polynomials with rational coefficients.
\newblock Mathematische Annalen. 1982 Dec;261(4):515-34.
\newblock Available from: \url{https://doi.org/10.1007/BF01457454}.

\end{thebibliography}
\end{document}